\documentclass[11pt]{article}
\usepackage{amsfonts,mathrsfs,amssymb,amsthm,mathptm}
\usepackage{amsmath,amscd,pslatex}
\usepackage{appendix}
\usepackage{multirow}
\usepackage{color}
\usepackage[all]{xy}
\usepackage{url}
\usepackage{exscale}
\usepackage{relsize}
\usepackage{bbm}
\usepackage{cite}

\begin{document}
\newtheorem{Def}{Definition}[section]
\newtheorem{Bsp}[Def]{Example}
\newtheorem{Prop}[Def]{Proposition}
\newtheorem{Theo}[Def]{Theorem}
\newtheorem{Lem}[Def]{Lemma}
\newtheorem{Koro}[Def]{Corollary}
\theoremstyle{definition}
\newtheorem{Rem}[Def]{Remark}

\newcommand{\add}{{\rm add}}
\newcommand{\con}{{\rm con}}
\newcommand{\gd}{{\rm gl.dim}}
\newcommand{\dm}{{\rm domdim}}
\newcommand{\tdim}{{\rm dim}}
\newcommand{\E}{{\rm E}}
\newcommand{\Mor}{{\rm Morph}}
\newcommand{\End}{{\rm End}}
\newcommand{\ind}{{\rm ind}}
\newcommand{\lcm}{{\rm lcm}}
\newcommand{\rsd}{{\rm res.dim}}
\newcommand{\rd} {{\rm rep.dim}}
\newcommand{\ol}{\overline}
\newcommand{\overpr}{$\hfill\square$}
\newcommand{\rad}{{\rm rad}}
\newcommand{\soc}{{\rm soc}}
\renewcommand{\top}{{\rm top}}
\newcommand{\stp}{{\mbox{\rm -stp}}}
\newcommand{\pd}{{\rm projdim}}
\newcommand{\id}{{\rm injdim}}
\newcommand{\fld}{{\rm flatdim}}
\newcommand{\fdd}{{\rm fdomdim}}
\newcommand{\Fac}{{\rm Fac}}
\newcommand{\Gen}{{\rm Gen}}
\newcommand{\fd} {{\rm findim}}
\newcommand{\Fd} {{\rm Findim}}
\newcommand{\Pf}[1]{{\mathscr P}^{<\infty}(#1)}
\newcommand{\DTr}{{\rm DTr}}
\newcommand{\Tr}{{\rm Tr}}
\newcommand{\cpx}[1]{#1^{\bullet}}
\newcommand{\D}[1]{{\mathscr D}(#1)}
\newcommand{\Dz}[1]{{\mathscr D}^+(#1)}
\newcommand{\Df}[1]{{\mathscr D}^-(#1)}
\newcommand{\Db}[1]{{\mathscr D}^b(#1)}
\newcommand{\C}[1]{{\mathscr C}(#1)}
\newcommand{\Cz}[1]{{\mathscr C}^+(#1)}
\newcommand{\Cf}[1]{{\mathscr C}^-(#1)}
\newcommand{\Cb}[1]{{\mathscr C}^b(#1)}
\newcommand{\Dc}[1]{{\mathscr D}^c(#1)}
\newcommand{\K}[1]{{\mathscr K}(#1)}
\newcommand{\Kz}[1]{{\mathscr K}^+(#1)}
\newcommand{\Kf}[1]{{\mathscr  K}^-(#1)}
\newcommand{\Kb}[1]{{\mathscr K}^b(#1)}

\newcommand{\modcat}{\ensuremath{\mbox{{\rm -mod}}}}
\newcommand{\Modcat}{\ensuremath{\mbox{{\rm -Mod}}}}
\newcommand{\stmodcat}[1]{#1\mbox{{\rm -{\underline{mod}}}}}
\newcommand{\Imodcat}[1]{#1\mbox{{\rm -Inj}}}
\newcommand{\imodcat}[1]{#1\mbox{{\rm -inj}}}
\newcommand{\Pmodcat}[1]{#1\mbox{{\rm -Proj}}}
\newcommand{\pmodcat}[1]{#1\mbox{{\rm -proj}}}
\newcommand{\GpModcat}[1]{#1\mbox{{\rm -GProj}}}
\newcommand{\Gpmodcat}[1]{#1\mbox{{\rm -Gproj}}}
\newcommand{\GiModcat}[1]{#1\mbox{{\rm -GInj}}}
\newcommand{\Gimodcat}[1]{#1\mbox{{\rm -Ginj}}}

\newcommand{\opp}{^{\rm op}}
\newcommand{\otimesL}{\otimes^{\rm\mathbb L}}
\newcommand{\rHom}{{\rm\mathbb R}{\rm Hom}\,}
\newcommand{\projdim}{\pd}
\newcommand{\Hom}{{\rm Hom}}
\newcommand{\Coker}{{\rm Coker}}
\newcommand{\Ker  }{{\rm Ker}}
\newcommand{\Cone }{{\rm Con}}
\newcommand{\Img  }{{\rm Im}}
\newcommand{\Ext}{{\rm Ext}}
\newcommand{\StHom}{{\rm \underline{Hom}}}

\newcommand{\gm}{{\rm _{\Gamma_M}}}
\newcommand{\gmr}{{\rm _{\Gamma_M^R}}}

\def\vez{\varepsilon}\def\bz{\bigoplus}  \def\sz {\oplus}
\def\epa{\xrightarrow} \def\inja{\hookrightarrow}
\newcommand{\lra}{\longrightarrow}
\newcommand{\llra}{\longleftarrow}
\newcommand{\lraf}[1]{\stackrel{#1}{\lra}}
\newcommand{\llaf}[1]{\stackrel{#1}{\llra}}
\newcommand{\ra}{\rightarrow}
\newcommand{\dk}{{\rm dim_{_{k}}}}

\newcommand{\colim}{{\rm colim\, }}
\newcommand{\limt}{{\rm lim\, }}
\newcommand{\Add}{{\rm Add }}
\newcommand{\Tor}{{\rm Tor}}
\newcommand{\Cogen}{{\rm Cogen}}
\newcommand{\mlt}{{\rm mlt}}
\newcommand{\Tria}{{\rm Tria}}
\newcommand{\tria}{{\rm tria}}

\newcommand{\mQ}{\mathsf{Q}}
\newcommand{\mP}{\mathsf{P}}
\newcommand{\mL}{\mathsf{L}}
\newcommand{\mR}{\mathsf{R}}
\newcommand{\mf}{\mathsf{F}}
\newcommand{\mg}{\mathsf{G}}
\newcommand{\mt}{\mathsf{T}}
\newcommand{\mh}{\mathsf{H}}
\newcommand{\mU}{\mathsf{U}}
\newcommand{\mz}{\mathsf{Z}}
\newcommand{\mc}{\mathsf{C}}
\newcommand{\mr}{\mathsf{r}}
\newcommand{\mq}{\mathsf{q}}
\newcommand{\mi}{\mathsf{i}}

\newcommand{\mD}{\mathsf{D}}

{\Large \bf
\begin{center}
Gorenstein projective modules over rings of Morita contexts
\end{center}}
	
\medskip
\centerline{\textbf{ Qianqian Guo and Changchang Xi$^*$}}
	
\renewcommand{\thefootnote}{\alph{footnote}}
\setcounter{footnote}{-1} \footnote{$^*$Corresponding author's
		Email: xicc@cnu.edu.cn; Fax: 0086 10 68903637.}
\renewcommand{\thefootnote}{\alph{footnote}}
\setcounter{footnote}{-1} 
\renewcommand{\thefootnote}{\alph{footnote}}
\setcounter{footnote}{-1} 
	
\begin{abstract} Under semi-weak and weak compatibility of bimodules, we establish sufficient and necessary conditions of Gorenstein-projective modules over rings of Morita contexts with one bimodule homomorphism zero. This generalises and extends results on triangular matrix Artin algebras and on Artin algebras of Morita contexts with two bimodule homomorphisms zero in the literature, where only sufficient conditions are given under a strong assumption of compatibility of bimodules. An application is provided to describe Gorenstein-projective modules over noncommutative tensor products arising from Morita contexts. Moreover, we work with Noether rings and modules instead of Artin algebras and modules.

\end{abstract}

{\footnotesize\tableofcontents\label{contents}}	
	
\section{Introduction}
For finitely generated modules over noetherian rings, Auslander and Bridge introduced Gorenstein-projective modules in \cite{AB69}, that is, modules of G-dimension zero, and this idea was generalized several decades later by Enochs, Jenda, and Torrecillas for arbitrary modules over arbitrary rings in
\cite{EJT93b, EJ95}. The notion of Gorenstein-projective modules plays nowadays a very important role in the so-called Gorenstein homological algebra which has significant applications in commutative algebra, algebraic geometry and other fields. It is fundamental, but also difficult, to describe all Gorenstein-projective modules over a given algebra or ring. Recently, there is a lot of interesting works done in this direction. For instance, in a series of articles \cite{LZ, XZ, zhangpu13}, Gorenstein-projective modules over the triangular matrix Artin algebra $\left(\begin{smallmatrix}A & N\\ 0 &B\end{smallmatrix}\right)$
were determined under some assumptions on the bimodule $_AN_B$. A natural generalization of triangular matrix algebras is the Morita context rings $\Lambda_{(0,0)}=\left(\begin{smallmatrix}A & N\\ M &B\end{smallmatrix}\right)$ with zero bimodule homomorphisms. Recall that a Morita context ring is generally the $2\times 2$ matrix ring $\Lambda_{(\phi,\psi)}:=\left(\begin{smallmatrix}A & N\\ M & B\end{smallmatrix}\right)$ associated to a Morita context $(A,B, {}_BM_A, {}_AN_B,\phi: M\otimes_AN\to B, \psi: N\otimes_BM\to A)$ with bimodule homomorphisms $(\phi, \psi)$, and that $\Lambda_{(\phi,\psi)}$-modules are presented by quadruples $(X,Y,f,g)$ (see Section \ref{sect2} for details). In \cite{GaoP17}, Gao and Psaroudakis give a set of concise  sufficient conditions for $\Lambda_{(0,0)}$-modules to be Gorenstein-projective. To achieve their results, they require some compatibility conditions on the bimodules $M$ and $N$, which were introduced in \cite{zhangpu13}.

In this note we consider Morita context rings $\Lambda_{\psi}:=\left(\begin{smallmatrix}A & N\\ M &B\end{smallmatrix}\right)_{(0,\psi)}$ with one bimodule homomorphism zero, and characterize their Gorenstein-projective modules.  To implement our characterization of Gorenstein-projective modules, we use weak versions of compatibility conditions (see Section \ref{sect2.2} for definition). With an additional assumption on a couple of very special $\Lambda_{\psi}$-modules related to ingredients of the given Morita context, we even can show that the weak compatibility conditions are necessary and sufficient for a $\Lambda_{\psi}$-module to be Gorenstein-projective. Moreover, we will work with Noether rings instead of Artin algebras. Our main results, Theorems  \ref{main-result-1} and \ref{main-result-3}, are summarised as follows.

\begin{Theo} \label{main-result-intr1} Let $A$ and $B$ be noetherian rings, and $(A, B, {}_BM_A, {}_AN_B,0, \psi)$ a Morita context with the bimodules $_BM_A$ and $_AN_B$ finitely generated as one-sided modules. Further, assume that $A$ is the trivial extension of a subring $\Lambda$ of $A$ by the image $I$ of $\psi$.

\smallskip
$(\rm{I})$ The following two sets of conditions are equivalent for the Morita context ring $\Lambda_{\psi}$.

$\quad (1)$ ${}_\Lambda N_B$, $_BM_\Lambda$ and $_{\Lambda}I_{\Lambda}$ are weakly compatible bimodules, $(_AN,0,0,0)$ and $(_AI,0,0,0)$ are semi-weakly compatible left $\Lambda_{\psi}$-modules; and $(M_A, 0,0,0)$ and $(I_A, 0,0,0)$ are semi-weakly compatible right $\Lambda_{\psi}$-modules.

$\quad (2)$ A $\Lambda_{\psi}$-module $(X,Y,f,g)$ is Gorenstein-projective if and only if

$\qquad (a)$  $_B\Coker(f)$ and $_{\Lambda}\Coker(g)$ are Gorenstein-projective, and	

$\qquad (b)$	$_B\Img(f)\simeq {}_BM\otimes_A\Coker(g)$, ${}_A\Img(g)/IX\simeq {}_AN\otimes_B\Coker(f)$, and
$_AIX\simeq {}_AI\otimes_A\Coker(g)$, where $\Coker(f)$ and $\Img(g)$ denote the cokernel of $f$ and the image of $g$, respectively.

\smallskip
$({\rm II})$ Suppose that $_\Lambda N_B$, $_BM_\Lambda$ and $_{\Lambda}I_{\Lambda}$ are weakly compatible. If a $\Lambda_{\psi}$-module $(X,Y,f,g)$ satisfies the above conditions $(a)$ and $(b)$ in $(2)$, then it is Gorenstein-projective.
\end{Theo}

The result (I) not only generalises and extends greatly the ones on triangular matrix algebras by Pu Zhang and his students, and on Morita context rings with two bimodule homomorphisms zero by Nan Gao and Chrysostomos Psaroudakis in \cite{zhangpu13, XZ} and \cite{GaoP17}, respectively, to a large class of Morita context algebras, but also can be applied to a class of noncommutative tensor products, see Corollary \ref{nonctp} for details. Notably, noncommutative tensor products generalize usual tensor products over
commutative rings, capture many known constructions in ring theory, and are useful in constructing reollements of derived module categories (see \cite{xc3, xc7}).

This note is structured as follows. In Section \ref{sect2} we recall (weakly) compatible bimodules, a complete Horseshoe lemma and basic facts on Morita context rings. In Section \ref{sect3} we prove the main result, Theorem \ref{main-result-intr1}, and then formulate it for the special Morita context rings $\Lambda_{(0,0)}$. In this case the resulting statement appears in a quite simple form, see Proposition \ref{zero-case}. Finally, in Section \ref{sect4} we apply our result to  noncommutative tensor products arising from Morita contexts with two bimodule homomorphisms zero. This provides in fact a corresponding result for Morita context rings $\Lambda_{(\phi,0)}$, as indicated by Corollary \ref{nonctp}.

\section{Preliminaries\label{sect2}}

In this section we recall basic definitions and facts for later proofs.

Let $A$ be a unitary (associative) ring. We denote by $A\Modcat$ (respectively, $A\modcat$) the category of all (respectively, finitely generated) left $A$-modules. As usual, $A\Pmodcat$ and $A\pmodcat$  are the full subcategories of $A\Modcat$ consisting of all projective modules and finitely generated projective $A$-modules, respectively. Similarly, we have the notation $A\Imodcat$ and $A\imodcat$ for the full subcategories of all injective $A$-modules and finitely generated injective $A$-modules, respectively.
For a full subcategory $\mathcal{X}$ of $A\Modcat$, we denote by $\mathscr{C}(\mathcal{X})$ the category of complexes over $\mathcal{X}$, and write $\mathscr{C}(A)$ for $\mathscr{C}(A\Modcat)$.

The composite of two homomorphisms $f: X\to Y$ and $g:Y\to Z$ will be denoted by $fg$ instead of $gf$. Thus the image of $x\in X$ under $f$ is written as $(x)f$ or $xf$, and the image of $f$ is denoted by $\Img(f)$.

Recall that a complex $\cpx{X}=(X^i,d_X^i)\in\mathscr{C}(A)$ is \emph{exact} if the cohomology group $H^i(\cpx{X})=0$ for all $i$; and \emph{totally exact} if it is exact and the complex $\Hom_A(\cpx{X},A)$ is exact. Let $X$ be an $A$-module. An exact complex $\cpx{P}\in \mathscr{C}(\Pmodcat{A})$ is called a \emph{complete projective resolution }of $X$ if $\Ker(d_P^0)= X$. By a \emph{total projective resolution} of $X$ we mean a totally exact, complete projective resolution of $X$.  Following \cite{EJ95}, the module $_AX$ is \emph{Gorenstein-projective} if $X$ has a total projective resolution. Dually, an $A$-module $Y$ is \emph{Gorenstein-injective} if there is a complete injective resolution $\cpx{I}\in\mathscr{C}(\Imodcat{A})$ such that $\Ker(d_I^0)=Y$ and $\Hom_A(E,\cpx{I})$
is exact for all $E\in \imodcat{A}$. In $A\modcat$, Gorenstein-projective modules are nothing else than modules of $G$-dimension $0$ in the sense of Auslander-Bridge \cite{AB69}.
We denote by $A\GpModcat$ (respectively, $\Gpmodcat{A})$ the category of all (respectively, finitely generated) Gorenstein-projective $A$-modules, and by $\GiModcat{A}$ (respectively, $\Gimodcat{A})$ the category of all (respectively, finitely generated) Gorenstein-injective $A$-modules. It is known that $\Gpmodcat{A}$ contains $\pmodcat{A}$ and is closed under direct summands, extensions and kernels of surjective homomorphisms (see \cite{holm}).

Since Gorenstein-projective modules involve complete projective resolutions, a kind of complete Horseshoe lemma will be useful. For the convenience of the reader, we state it here for module categories and still refer to Horseshoe Lemma. For other versions, see \cite{holm}, \cite{zhangpu13} and \cite{GaoP17}.

\begin{Lem}[Horseshoe Lemma] \label{horseshoe} Given a short exact sequence $0\to U\to W\to V\to 0$ of $A$-modules and two exact complexes $\cpx{X}=(X^i,d_X^i)$ and $\cpx{Y}=(Y^i,d_Y^i)$ of $A$-modules with $\Ker(d_X^0)=U$ and $\Ker(d_Y^0)=V$, if $\Ext^1_A(\Ker(d^i_Y),X^i)=0$ for all $i\ge 0$ and if $\Ext^1_A(Y^{-i},\Img(d_X^{-i}))=0$ for all $i\ge 1$, then there is an exact complex $\cpx{Z}=(Z^i,d_Z^i)$ and an exact sequence of complexes
$$ 0\lra \cpx{X}\lra\cpx{Z}\lra \cpx{Y}\lra 0,$$ where $Z^i=X^i\oplus Y^i$, $d_Z^i=\big(\begin{smallmatrix}d^i_X & 0 \\ \rho^i & d^i_Y\end{smallmatrix}\big),$ and $\rho^i: Y^i\ra X^{i+1}$ is a homomorphism of $A$-modules, such that the induced exact sequence $0\ra \Ker(d^0_X)\ra \Ker(d_Z^0)\ra \Ker(d_Y^0)\ra 0$ coincides with the given short exact sequence.

Further, if $X^i=X^{i+1}$, $Y^i=Y^{i+1}$,$d^i_X=d_X^{i+1}$ and $d^i_Y=d_Y^{i+1}$ for all $i$, then $Z^i=Z^{i+1}$  and $d_Z^i=d_Z^{i+1}$  for all $i$.
\end{Lem}

The following easy lemma is often used, its proof is left to the reader.

\begin{Lem} \label{ker-seq-sur} $(1) $ If $0\ra \cpx{C} \lraf{\cpx{c}} \cpx{E} \lraf{\cpx{d}}\cpx{G} \ra 0$ is an exact sequence of complexes of $A$-modules, then there is an induced  exact sequence
$  0\ra  \Ker(d_C^i) \lraf{\bar{c}^i}  \Ker(d_E^i) \lraf{\bar{d}^i} \Ker(d_G^i)  \; $
in $A\Modcat \,$ for any $i \in \mathbb{Z}$. In particular, if $\cpx{C}$ and $\cpx{E}$ are exact, then $\bar{d}^i$ is surjective.

$(2)$ Given a $3$-dimensional diagram of $A$-modules with the squares consisting of solid arrows commutative
$$\xymatrix@C=0.5cm@R=0.3cm{
	& K' \ar[dl] \ar[rr] \ar'[d][dd]
	&  & C' \ar'[d][dd]  \ar[rr]\ar[dl] & & A' \ar[dd]  \ar[dl]   \\
	L' \ar[rr]^(0.77){l'}\ar@{-->}[dd]
	&  &  D'  \ar[dd] \ar[rr]^(0.77){d'} & & B' \ar[dd] \\
	&  K \ar@{-->}[dl]\ar'[r]^(0.5){k}[rr]
	&  &    C  \ar[dl]\ar'[r]^(0.5){c}[rr] & &    A\ar[dl]      \\
	L \ar[rr]^{l }
	&  & D \ar[rr]^{d} & &  B  }$$
If $l'd'=0= kc$ and $l$ is the kernel of $d$, then the dashed arrows exist and every new square commutes.
\end{Lem}

The next lemma is well known.

\begin{Lem}\label{noether} Let $R$ be a left noetherian ring, and $M, N$ be finitely generated $R$-modules.

$(1)$ Every surjective homomorphism $f: {}_RM\to {}_RM$ is an automorphism.

$(2)$ If $M\simeq N$ and $f: M\to N$ is a surjective homomorphism of $R$-modules, then $f$ is an isomorphism.
\end{Lem}

Finally, we recall the notion of approximations. Let $\mathcal{D}$  be a full additive subcategory of an additive category $\mathcal{C}$ and $X$ an object in $\mathcal{C}$. A morphism $f: X\to D$ in $\mathcal{C}$ is called a \emph{left $\mathcal{D}$-approximation} of $X$ if $D\in \mathcal{D}$ and $\Hom_{\mathcal{C}}(f,D'): \Hom_{\mathcal{C}}(D,D')\to \Hom_{\mathcal{C}}(X,D')$ is surjective for any object $D'\in\mathcal{D}$. Dually, a morphism $f: D\to X$ in $\mathcal{C}$ is called a \emph{right $\mathcal{D}$-approximation} of $X$ if $D\in \mathcal{D}$ and $\Hom_{\mathcal{C}}(D',f): \Hom_{\mathcal{C}}(D',D)\to \Hom_{\mathcal{C}}(D', X)$ is surjective for any object $D'\in\mathcal{D}$. Remark that left and right approximations are termed as preenvelopes and precovers in ring theory, respectively.

\subsection{Morita context rings and their modules\label{sect2.1}}

Morita context rings stem from description of Morita equivalences of rings (see \cite{morita, bass}) and appear in many situations (for example, see \cite{GP}). There is a large variety of literature on Morita contexts, duality and equivalences (for instance, \cite{morita, bass, buchweitz, MR, Green, KT, GaoP17}). We briefly recall  Morita context rings and their modules.

Let $A$ and $B$ be unitary rings, and let $_AN_B$ be an $A$-$B$-bimodule, $_BM_A$ a $B$-$A$-bimodule, $\phi: M\otimes_{A}N \ra B$ a homomorphism of $B$-$B$-bimodules and $\psi : N\otimes_{B}M \ra A$ a homomorphism of $A$-$A$-bimodules. Recall that the sextuple $(A,B, M, N, \phi, \psi)$ is a \emph{Morita context} (see \cite{morita}) if
the following two diagrams are commutative
$$\xymatrix{ N\otimes_BM\otimes_AN \ar[r]^(0.6){1_N\otimes \phi}\ar[d]^{\psi\otimes 1_N} & N\otimes_BB\ar[d]^{\mlt} \\ A\otimes_AN \ar[r]^{\mlt} & N, } \qquad
\xymatrix{ M\otimes_AN\otimes_BM \ar[r]^(0.6){\phi\otimes 1_M }\ar[d]^{1_M\otimes \psi} & B\otimes_BM\ar[d]^{\mlt} \\ M\otimes_AA \ar[r]^{\mlt} & M }$$where $\mlt$ stands for the multiplication map universally.

Let $I:=\Img(\psi) $ and $J:=\Img(\phi)$. Associated with a Morita context $(A,B, M, N, \phi, \psi)$, there is defined a \emph{Morita context ring} (see \cite{morita, bass}), denoted by $\Lambda_{(\phi,\psi)}$, which has the underlying abelian group of the matrix form
with the multiplication induced by $\phi$ and $\psi$:
$$\Lambda_{(\phi,\psi)}:=         \left(\begin{matrix}
           A & N \\
           M & B \\
         \end{matrix}\right)= \left\{\left(\begin{matrix}
           a & n \\
           m & b \\
         \end{matrix}\right) \ | \ a\in A, b \in B, n\in N, m\in M\right\},$$  $$\left(\begin{matrix}
           a & n \\
           m & b \\
         \end{matrix}\right)\left(\begin{matrix}
           a' & n' \\
           m' & b' \\
         \end{matrix}\right)=\left(\begin{matrix}
           aa'+(n\otimes m')\psi & an'+nb' \\
           ma'+bm' & (m\otimes n')\phi +bb' \\
         \end{matrix}\right).                           $$
We write $\Lambda_{\psi}=\Lambda_{(0,\psi)}$ and $\Lambda_{\phi}=\Lambda_{(\phi,0)}$.

The description of modules over $\Lambda_{(\phi,\psi)}$ are well known, for instance, see \cite{Green, KT}. We recall it here very briefly. Every $\Lambda_{(\phi,\psi)}$-module is determined by a quadruple $(X,Y,f,g)$, where $X$ and $Y$ are modules over $A$ and $B$, respectively, $f\in\Hom_B(M\otimes_AX,Y)$ and $g\in \Hom_A(N\otimes_BY,X)$ such that the diagrams are commutative
\begin{equation}
\xymatrix{
  N\otimes_B M\otimes_A X \ar[d]_{\psi\otimes 1_X} \ar[r]^{ \ \ \ \ \ 1_N\otimes f} &  N\otimes_BY \ar[d]^{g}     \\
  A\otimes_AX    \ar[r]^{\simeq} & X ,                 } \ \ \ \ \ \ \  \ \ \  \xymatrix{
  M\otimes_A N\otimes_B Y \ar[d]_{\phi\otimes 1_Y} \ar[r]^{ \ \ \ \ 1_M\otimes g} &  M\otimes_AX \ar[d]^{f}     \\
  B\otimes_BY    \ar[r]^{\simeq} & Y                  }
\end{equation}
where the two isomorphisms are the multiplication maps.

If $(_AX, {}_BY, f, g)$ is a $\Lambda_{(\phi,\psi)}$-module, then $I\Coker(g)=0$ and $J\Coker(f)=0$. This follows from (1) since $IX\subseteq \Img(g)$ and $JY\subseteq \Img(f)$. Thus $\Coker(g)$ is an $A/I$-module and $\Coker(f)$ is a $B/J$-module.
Since $\Hom_B(M\otimes_AX,Y)\simeq \Hom_A(X,\Hom_B(M,Y))$, we denote by $\tilde{f}: X\to \Hom_B(M,Y)$ the image of $f$ under this adjunction. Similarly, we define $\tilde{g}: Y\to \Hom_A(N,X)$.

A homomorphism from a $\Lambda_{(\phi,\psi)}$-module $(X,Y,f,g)$ to another $\Lambda_{(\phi,\psi)}$-module $(X',Y',f',g')$ is a pair $(\alpha, \beta)$ with $\alpha\in \Hom_A(X,X')$ and $\beta\in \Hom_B(Y,Y')$ such that the two diagrams are commutative
\begin{equation}\label{morphismdiagrams}
\xymatrix{
  M\otimes_A X \ar[d]_{1_{M}\otimes \alpha} \ar[r]^{\ \ \ f} & Y \ar[d]^{\beta}     \\
  M\otimes_AX'    \ar[r]^{\ \ \ f'} & \; Y'}\qquad \mbox{and}  \qquad \xymatrix{
  N\otimes_B Y \ar[d]_{1_{N}\otimes \beta} \ar[r]^{\ \ \  g} &  X \ar[d]^{\alpha}     \\
  N\otimes_BY'    \ar[r]^{\ \ \ g'} & X'                  }
\end{equation}
Clearly, for a homomorphism $(\alpha, \beta): (X,Y,f,g)\to (X',Y',f',g')$ of $\Lambda_{(\phi,\psi)}$-modules, its kernel  $\Ker(\alpha,\beta)$ is $(\Ker(\alpha),\Ker(\beta), h, j)$, where $h$  and $j$ are uniquely given by the commutative diagrams:
\begin{equation}
\label{diagramskernels}
 \xymatrix{
    M\otimes_A\Ker(\alpha) \ar@{-->}[d]_{h} \ar[r]^{ \ \ 1_{M}\otimes i_X} & M\otimes_{A}X \ar[d]_{f}
  \ar[r]^{ 1_{M}\otimes \alpha} & M\otimes_A X' \ar[d]^{f'} \\
   \Ker(\beta) \ \ar@{^{(}->}[r]^{i_Y} & Y \ar[r]^{\beta} & \; Y', } \ \ \ \ \  \xymatrix{
    N\otimes_B\Ker(\beta) \ar@{-->}[d]_{j} \ar[r]^{ \ \ 1_N\otimes i_Y} & N\otimes_BY \ar[d]_{g}
  \ar[r]^{1_{N}\otimes\beta} & N\otimes_BY' \ar[d]^{g'} \\
   \Ker(\alpha) \ \ar@{^{(}->}[r]^{i_X} & X \ar[r]^{\alpha} & X'   }
\end{equation}
where $i_X: \Ker(\alpha)\to X$ and $i_Y:\Ker(\beta)\to Y$ are the inclusions. Dually, one describes the cokernel of $(\alpha,\beta)$.

Let $0 \ra (X_1,Y_1,f_1,g_1) \lraf{(\alpha_1,\beta_1)} (X_2,Y_2,f_2,g_2) \lraf{(\alpha_2,\beta_2)} (X_3,Y_3,f_3,g_3) \ra 0$ be a sequence of $\Lambda_{(\phi,\psi)}$-modules. This sequence is exact if and only if the induced sequences $0\ra X_1\lraf{\alpha_1} X_2\lraf{\alpha_2} X_3\ra 0$ and $0\ra Y_1\lraf{\beta_1} Y_2\lraf{\beta_2} Y_3\ra 0$ are exact in $A\Modcat$ and $B\Modcat$, respectively.

Given two $\Lambda_{(\phi,\psi)}$-modules $(X,Y,f,g)$ and $(X',Y',f',g')$, their direct sum is given by $(X\oplus X', Y\oplus Y', f\oplus f', g\oplus g')$, where
$ f \oplus f' = \left(\begin{smallmatrix}
f & 0 \\
0 & f' \\
\end{smallmatrix}\right): M\otimes_AX\, \oplus \,M\otimes_A X'\ra Y\oplus Y'$ is defined to be the diagonal homomorphism of $B$-modules.

For $X\in A\Modcat$ and $Y\in B\Modcat$, we denote by $\Psi_X$ and $\Phi_Y$ the composites of the maps, respectively:
\[
\xymatrix@C=0.5cm{
  N\otimes_BM\otimes_AX \ar[rr]^{ \ \ \ \ \psi\otimes 1_X}  \ar @/^1.5pc/[rrrr]^{{\Psi_X}} && A\otimes_AX  \ar[rr]^{\simeq} && X,  } \ \ \ \ \xymatrix@C=0.5cm{
  M\otimes_AN\otimes_BY \ar[rr]^{ \ \ \ \ \phi\otimes 1_Y}  \ar @/^1.5pc/[rrrr]^{{\Phi_Y}} && B\otimes_BY  \ar[rr]^{\simeq} && Y.  }
\]

The bimodules $_BM_A$ and $_AN_B$ define two natural transformations $\zeta$ and $\xi$  between tensor functors and hom-functors:
$$M\otimes_A- \lraf{\zeta} \Hom_A(N,-): A\Modcat\lra B\Modcat, $$ $$ \zeta_X: M\otimes_AX\lra \Hom_A(N,X), \quad m\otimes x \mapsto [n \mapsto (n\otimes m)\psi\, x],   $$
$$N\otimes_B- \lraf{\xi} \Hom_B(M,-): B\Modcat\lra A\Modcat, $$ $$\xi_Y:N\otimes_BY \ra \Hom_B(M,Y), \quad n\otimes y\mapsto [ m \mapsto (m\otimes n)\phi\, y].$$
Following \cite{GP}, we define functors related to Morita context rings as follows.
$$\mt_A: A\Modcat\lra \Lambda_{(\phi,\psi)}\Modcat, \quad {}_AX\mapsto \mt_{A}(X)=(X,M\otimes_AX,1_{M\otimes X},\Psi_X),$$
$$\mh_A: A\Modcat\lra \Lambda_{(\phi,\psi)}\Modcat, \quad {}_AX\mapsto \mh_A(X)=(X,\Hom_{A}(N,X),\zeta_{X},\delta_X), $$
$$\mt_B: B\Modcat\lra \Lambda_{(\phi,\psi)}\Modcat, \quad {}_BY\mapsto \mt_{B}(Y)=(N\otimes_BY,Y,\Phi_Y,1_{N\otimes Y}),$$
$$ \mh_B: B\Modcat\lra \Lambda_{(\phi,\psi)}\Modcat, \quad {}_BY\mapsto \mh_B(Y)=(\Hom_{B}(M,Y),Y,\delta_{Y},\xi_Y),$$ where $\delta_X: N\otimes_B\Hom_A(N,X)\ra X$ and $\delta_Y: M\otimes_A\Hom_B(M,Y)\ra Y$ are evaluation maps. Note that for an $A$-module $X$ with $IX=0$ and a $B$-module $Y$ with  $JY=0$, we can get naturally  $\Lambda_{(\phi,\psi)}$-modules $(X, 0,0,0)$ and $(0,Y,0,0),$ respectively. This gives rise to the functors
$$\mz_{A/I}: (A/I)\Modcat\lra \Lambda_{(\phi,\psi)}\Modcat, \quad {}_{A/I}U\mapsto \mz_{A/I}(U)=(_AU,0,0,0),$$
$$\mz_{B/J}: (B/J)\Modcat\lra \Lambda_{(\phi,\psi)}\Modcat, \quad {}_{B/J}V\mapsto \mz_{B/J}(V)=(0,{}_BV,0,0).$$
The actions of all functors defined above on morphisms are defined naturally. The relation among these functors are given by the following lemma. Since our proofs use only adjoint pairs of functors in recollements, we will not recall the definition of recollements of abelian categories, and just refer the reader to \cite{fp04, psaroudakis} for more details.

\begin{Lem}{\rm \cite{GP}} \label{mr-recol} There are two recollements of module categories:
\[
\xymatrix@C=0.5cm{
B/J\Modcat \ar[rrr]^{\mz_{B/J} \ \ } &&& \Lambda_{(\phi,\psi)}\Modcat \ar[rrr]^{\mU_A } \ar @/_1.5pc/[lll]_{\mathsf{Q}_B}  \ar
 @/^1.5pc/[lll]_{\mathsf{P}_B} &&& A\Modcat,
\ar @/_1.5pc/[lll]_{\mt_A} \ar
 @/^1.5pc/[lll]_{\mh_A}
 }
\]
\[
\xymatrix@C=0.5cm{
A/I\Modcat \ar[rrr]^{\mz_{A/I} \ \ } &&& \Lambda_{(\phi,\psi)}\Modcat \ar[rrr]^{\mU_B } \ar @/_1.5pc/[lll]_{\mathsf{Q}_A}  \ar
 @/^1.5pc/[lll]_{\mathsf{P}_A} &&& B\Modcat
\ar @/_1.5pc/[lll]_{\mt_B} \ar
 @/^1.5pc/[lll]_{\mh_B}
 }\]
where $\mathsf{U}_A$ and $\mathsf{U}_B$ are the canonical projections to $A\Modcat$ and $B\Modcat$, respectively, where $\mathsf{Q}_A= ((A/I)\otimes_A-) \mathsf{U_A}$ and $\mathsf{Q}_B= ((B/J)\otimes_B-) \mathsf{U_B}$, and where $\mathsf{P}_A$ and $\mathsf{P}_B$ are defined on objects $(X,Y,f,g)$ by taking kernels of $\tilde{f}$ and $\tilde{g}$, respectively.
\smallskip
\end{Lem}

Suppose that $A,B$ are noetherian rings and $_BM_A,{}_AN_B$ are bimodules such that they are finitely generated as one-sided modules. Then it is known that $\Lambda_{(\phi,\psi)}$ is a noetherian ring (see, for example,\cite[Proposition 1.7, p.12]{MR}). For a noetherian ring, its identity has a complete decomposition of orthogonal primitive idempotent elements (see [Proposition 10.14, p.128] \cite{af74}). Thus the description of indecomposable projective modules over the Artin algebra $\Lambda_{(\phi,\psi)}$ in \cite[Proposition 3.1]{GP} extends to the one over the noetherian ring $\Lambda_{(\phi,\psi)}$.

\begin{Lem}\label{prop:projmod} Suppose that $A,B$ are noetherian rings and $_BM_A,{}_AN_B$ are bimodules such that they are finitely generated as one-sided modules.

$(1)$ {\rm \cite[Proposition 3.1]{GP}} An indecomposable $\Lambda_{(\phi,\psi)}$-module is projective if and only if it is given by
$   \mt_{A}(P)=(P, M\otimes_{A}P, id_{M\otimes_AP}, \Psi_P)$, or
$   \mt_B(Q)=({N\otimes_BQ}, Q, \Phi_{Q}, id_{N\otimes_{B}Q})$, where $P$ and $Q$ are indecomposable projective modules over $A$ and $B$, respectively.

$(2)$ {\rm \cite[Corollary 2.2]{mueller}} An indecomposable $\Lambda_{(\phi,\psi)}$-module is injective if and only if it is of the form
$   \mh_A(U)=(U,\Hom_{A}(N,U),$ $\zeta_{U},\delta_U)$, or
$   \mh_B(V)=(\Hom_{B}(M,V),V,\delta_{V},\xi_V)$, where $U$ and $V$ are indecomposable injective modules over $A$ and $B$, respectively. 
\end{Lem}

\subsection{Weakly and semi-weakly compatible modules\label{sect2.2}}

Compatible modules were defined in \cite{zhangpu13} to describe a class of Gorenstein-projective modules for triangular matrix Artin algebras which are of course special Morita context rings. They were further pursued in \cite{GaoP17} for Morita context Artin algebras $\Lambda_{(0,0)}$. We will use weakly and semi-weakly compatible modules for constructing Gorenstein-projective modules over noetherian rings $\Lambda_{\psi}$ that are more general than $\Lambda_{(0,0)}$.

Let $A$ and $B$ be unitary rings. First, we recall the definition of (weakly) compatible bimodules.

\begin{Def} Let $_AN_B$ be a bimodule.

$(1)$ $_A N_B$ is compatible {\rm \cite[Definition 1.1]{zhangpu13}} if

$\quad (C1)$ $\Hom_A(\cpx{P},N)$ is exact for all totally exact complex $\cpx{P}\in \mathscr{C}(A\mbox{-proj})$, and

$\quad (C2)$	$N\otimes_B\cpx{Q}$ is exact for all exact complex $\cpx{Q} \in \mathscr{C}(B\mbox{-proj})$.

$(2)$ $_AN_B$ is \emph{weakly compatible} {\rm \cite[Definition 4.1]{HZ}	} if it satisfies
$(C1)$ and

$\quad (C3)$	$N\otimes_B\cpx{Q}$ is exact for all totally exact complex $\cpx{Q} \in \mathscr{C}(B\mbox{-proj})$.
\end{Def}

Note that weakly compatible bimodules require exactness only for totally exact complexes $\cpx{Q} \in \mathscr{C}(B\mbox{-proj})$ in $(C2)$, and that the notion of weakly compatible bimodules is a proper generalisation of the one of compatible bimodules (see \cite[Example 4.3]{HZ})).

\begin{Def} \label{weak}
$(1)$ A right $B$-module $Y_B$ is \emph{semi-weakly compatible} if $Y$ satisfies only $(C3)$, that is, $Y\otimes_B\cpx{Q}$ is exact for all totally exact complex $\cpx{Q} \in \mathscr{C}(B\mbox{-proj})$.

$(2)$ A left $A$-module $_AX$ is \emph{semi-weakly compatible} if it satisfies $(C1)$, that is, $\Hom_A(\cpx{P}, X)$ is exact for all totally exact complex $\cpx{P}\in \mathscr{C}(A\mbox{-proj})$.
\end{Def}

Transparently, a bimodule $_AN_B$ is weakly compatible if and only if $_AN$ and $N_B$ are semi-weakly compatible. For an Artin algebra $A$, there is a duality $D$: $A$-mod $\to$ $A^{\opp}$-mod. Thus a left $A$-module $X \in A$-mod is semi-weakly compatible if and only if the right $A$-module $D(X)_A$ is semi-weakly compatible.

\begin{Lem}\label{not-cpa} Let $A, B$ and $C$ be rings.

$(1)$ An $A$-module $_AX$ is semi-weakly compatible if and only if the right $A^{\opp}$-module $X_{A^{\opp}}$ is semi-weakly compatible.

$(2)$ Let $_AN_B$ be an $A$-$B$-bimodule.

$\quad (i)$ If the modules $_AN$ and $N_B$ are of finite injective dimension, then $_AN_B$ is a weakly compatible $A$-$B$-bimodule.

$\quad (ii)$ If there is an exact complex $\cpx{P}=(P^i,d_P^i)\in \mathscr{C}(\pmodcat{B})$
and an integer $i$ such that $\Tor_{1}^{B}(N,\Ker(d_{P}^i))\neq 0$, then $_AN_B$ is not compatible.

$(3)$ If bimodules $_AX_B$ and $_{C^{op}}Y_{B^{op}}$ are compatible, then the bimodule $_AX \otimes_{B}Y_C$ is weakly compatible.
\end{Lem}

{\bf Proof.} (1)  For a finitely generated projective module $_AP$, there is an isomorphism $\Hom_A(P, {}_AU)\otimes_BV\simeq \Hom_A(P, {}_AU\otimes_BV)$ as abelian groups for any bimodule $_AU_B$ and $B$-module $_BV$. If $\cpx{P}\in \mathscr{C}(\pmodcat{A} )$ is totally exact, then so is $\Hom_A(\cpx{P},A)\in \mathscr{C}(\pmodcat{A^{\opp}})$. Hence $_AX$ is semi-weakly compatible if and only if $\Hom_A(\cpx{P},X)$ is exact if and only if $\Hom_A(\cpx{P},A)\otimes_AX $ is exact if and only if $X\otimes_{A^{\opp}}\Hom_A(\cpx{P},A)$ is exact for all totally exact complex in $\mathscr{C}(\pmodcat{A})$. Note that $\Hom_A(-,A)$ is a duality between $A\pmodcat$ and $A^{\opp}\pmodcat$. Thus (1) holds.

(2)  $(i)$ If $_AN$ is injective, then it is semi-weakly compatible. Suppose that $_AN$ is of injective dimension $n$ and $0\ra N\to I\to X\to 0$ be an exact sequence with $I$ injective and $_AX$ of injective dimension $ n-1$. Then $_AX$ is semi-weakly compatible by induction. For a totally exact complex $\cpx{P}\in  \mathscr{C}(A\mbox{-proj})$, since $P^i$ is projective, we have an exact sequence of complexes:
$ 0\to \Hom_A(\cpx{P}, N)\to \Hom_A(\cpx{P},I_0)\to \Hom_A(\cpx{P},X)\to 0$ with both $\Hom_A(\cpx{P},I_0)$ and $ \Hom_A(\cpx{P},X)$ being exact. Thus $\Hom_A(\cpx{P}, N)$ is exact,  and therefore $_AN$ is semi-weakly compatible. Now, suppose that $N_B$ has finite injective dimension. This means that the left $B^{\opp}$-module $_{B^{\opp}}N$ has finite injective dimension. Thus $_{B^{\opp}}N$ is semi-weakly compatible. By (1), the right $B$-module $N_B$ is semi-weakly compatible. Hence $(i)$ follows.

$(ii)$ The exact sequence $0 \ra \Ker(d_{P}^{i-1}) \ra P^{i-1} \ra  \Ker(d_{P}^{i}) \ra 0 $ shows that
$\Tor_{1}^{B}(N,\Ker(d_{P}^i))\neq 0$ and $N\otimes_{B} \Ker(d_{P}^{i-1}) \ra N\otimes_{B}P^{i-1} $ is not an injective homomorphism. This means that the complex $N\otimes_B \cpx{P}$ is not exact in degree $i-1$, and therefore $_AN_B$ is not a compatible $A$-$B$-bimodule.

(3) Suppose that $\cpx{P}\in  \mathscr{C}(A\mbox{-proj})$  is a totally exact complex. Then each $P^i$ is finitely generated projective module, and there is an isomorphism $\Hom_{A}(P^i,X\otimes_{B}Y)\simeq \Hom_{A}(P^i,X)\otimes_{B}Y$.
This yields an isomorphism of complexes: $\Hom_A(\cpx{P},X\otimes_BY)\simeq \Hom_A(\cpx{P},X)\otimes_BY$.
Since $_AX$ is semi-weakly compatible, the complex $\Hom_{A}(\cpx{P},X)$ is exact. Thus $_{C^{\opp}}Y\otimes_{B^{\opp}}\Hom_A(\cpx{P},X)$ is exact by the compatibility of $_{C^{op}}Y_{B^{\opp}}$, that is, $\Hom_{A}(\cpx{P},X)\otimes_BY$ is exact. Hence $_AX\otimes_{B}Y$ is a semi-weakly compatible left $A$-module.

Let $\cpx{Q}\in \mathscr{C}(C\mbox{-proj})$ be  a totally exact complex. Since $_{C^{op}}Y_{B^{op}}$ are compatible, $_{C^{op}}Y$ is semi-weakly compatible. By (1), $Y_C$ is semi-weakly compatible. Thus $Y\otimes_C\cpx{Q}$ is exact. Since $_AX_B$ is compatible, the complex $X\otimes_B(Y\otimes_C\cpx{Q})$ is exact. Thus $X\otimes_BY_C$ is a semi-weakly compatible right $C$-module.
$\square$

\section{Proof of the main result \label{sect3}}

In the rest of this paper we assume that all rings considered are noetherian, that is both left and right noetherian, and all modules are finitely generated.

Let $(A,B, {}_BM_A, {}_AN_B,\phi,\psi)$ be a Morita context with $\phi=0$. We  consider the Morita context ring $\Lambda_{\psi} :=\left(\begin{smallmatrix}A & N\\ M &B\end{smallmatrix}\right)_{(0,\psi)}$. Let $I:=\Img(\psi)$. Then $IN=MI=0$ and $I^2=0$.
Assume further that $\Lambda$ is a subring of $A$ with the same identity and $A$ is the \emph{trivial extension} of $\Lambda$ by $I$, that is $A=\Lambda \ltimes I$ with the multiplication
$$(\lambda, x)(\lambda',x')=(\lambda\lambda',\lambda x'+ x\lambda'), \quad \lambda,\lambda'\in\Lambda,\,  x, x'\in I. $$
Thus $\Lambda \simeq A/I$. Let $\pi: A\ra \Lambda$ be the canonical surjection. The restriction of $\pi$ on $\Lambda$ is the identity $id_{\Lambda}$. Clearly, $I$ is an ideal of $A$ with $I^2=0$.

Every $A$-module restricts to a $\Lambda$-module via the inclusion of $\Lambda$ into $A$. Conversely, every $\Lambda$-module $X$ induces an $A$-module $A\otimes_{\Lambda}X= X\oplus I\otimes_{\Lambda}X$, and restricts to $A$-module via $\pi$, that is, by defining $IX=0$.

For a  $\Lambda_{\psi}$-module $(X,Y,f,g)$, let $\lambda_X: X\to \Coker(g)$ and $\mu_Y: Y\to \Coker(f)$ be the canonical projections.

\subsection{Sufficient conditions for Gorenstein-projective modules\label{sect3.1}} 

We first prove the following lemma.

\begin{Lem} \label{module-A}
If $(_AX, {}_BY, f, g)\in \Lambda_{\psi}\modcat$, then

$(1)$  $I \Coker(g)=0$ and $I\, \Img(g)=0$.

$(2)$  There is a unique $B$-module homomorphism $\eta_{Y}:M\otimes_{A} \Coker(g) \ra Y$ such that $f=(1_M\otimes \lambda_{X})\eta_Y$. Thus $\Img(f)=\Img(\eta_Y)$ and $\Coker(f)=\Coker(\eta_Y)$.

$(3)$ Let $p_X: X\ra X/IX$ be the canonical projection. Then there is a unique homomorphism $\theta_{X}:N\otimes_{B} \Coker f\ra X/IX $ of $A$-modules such that $g p_{X}=(1_N\otimes \mu_{Y})\theta_X$. Thus $\Img(\theta_{X}) = \Img( g p_X)$.

$(4)$ Let $\emph{mlt}_{X}: I\otimes_{A}X \ra X$ be the multiplication map. Then there is a unique homomorphism $m_{X}:I\otimes_{A} \Coker g \ra X$ of $A$-modules such that $\Img(m_{X})=IX$
 and $\mlt_{X}=(1_I\otimes \lambda_{X})m_X$.
\end{Lem}

{\bf Proof.} (1) Clearly, $I\Coker(g)=0$ holds for any Morita context rings and their modules. It follows from $IN=0$ that $I(N\otimes_BY)=0$. Thus
$I\Img(g)=I(N\otimes_BY))\lambda_X=(I(N\otimes_BY))\lambda_X=0.$

(2) There is the exact commutative diagram of $B$-modules
$$ \xymatrix{_BM\otimes_AN\otimes_BY\ar[r]^-{1_M\otimes g}\ar[d]^-{\phi\otimes 1_Y} & {}M\otimes_AX \ar[r]^-{1_M\otimes\lambda_X}\ar[d]^-{f} & M\otimes_A\Coker(g)\ar[r]\ar@{-->}[dl]^-{\exists ! \; \eta_Y} & 0\\
 B\otimes_BY\ar[r]^-{\mlt}& Y& &}
$$Since $\phi=0$, there is a unique homomorphism $\eta_Y$ of $B$-modules such that $f=(1_M\otimes\lambda_X)\eta_Y$. This implies that $\Img(f)=\Img(\eta_Y)$ because $1_M\otimes\lambda_X$ is surjective. Thus it follows from Lemma \ref{noether} that $\eta_Y$ is injective if and only if $M\otimes_A\Coker(g)\simeq \Img(f)$ as $B$-modules.

(3) There is the following commutative diagram of $A$-modules
$$\xymatrix@C=0.8cm@R=0.8cm{
  N\otimes_B M\otimes_A X \ar[d]_{\psi\otimes 1_X} \ar[r]^{ \ \ \ \ \ 1_N\otimes f} &  N\otimes_BY \ar[d]^{g} \ar[r]^-{1_N\otimes \mu_Y} &  N\otimes_B \Coker(f) \ar@{-->}[d]^{\exists ! \; \theta_X} \ar[r] & 0    \\
  IX \ar[d]   \ar@{^(->}[r] & X  \ar[r]^{p_X}\ar[d]^-{\lambda_X}  & X/IX  \ar[r]\ar[d]^-{\lambda'_X} &     0        \\
 0 \ar[r]& \Coker(g) \ar@{=}[r] & \Coker(g)
   }$$
It follows now from Lemma \ref{noether} that $\theta_X$ is injective if and only if $N\otimes_B\Coker(f)\simeq \Img(g)/IX.$

(4) Consider the commutative diagram of $A$-modules
$$  \xymatrix{
I\otimes_{A}\Img(g) \ar[r]& I\otimes_{A}X \ar[d]^{\mlt_{X}}\ar[r]^-{1_I\otimes \lambda_X} & I\otimes_{A}\Coker(g) \ar@{-->}[dl]^{\exists ! \; m_X}\ar[r] &0\\
& X
}$$
Thanks to $I\,\Img(g)=0$, there is a unique homomorphism $m_{X}$ of $A$-modules such that $\mlt_X= (1_I\otimes\lambda_X) m_X.$
Thus $m_X$ is injective if and only if $IX\simeq I\otimes_A\Coker(g)$ as $A$-modules by Lemma \ref{noether}. $\square$

\begin{Lem} \label{b-condition}
If $(X,Y,f,g) \in \Lambda_{\psi} \modcat$ satisfies
$_AN\otimes_B\Coker(f)\simeq {}_A\Img(g)/IX$, and
$_AIX\simeq {}_AI\otimes_A\Coker(g)$, then $\Img(g)$ is the pushout of $1_N\otimes_B \eta_Y$ and $\psi\otimes_{\Lambda} 1_{\Coker(g)}$, where $\eta_Y$ is given in Lemma  $\ref{module-A}(2)$.
\end{Lem}

{\bf Proof.} Put $U:=\Coker(g)$, $V:=\Coker(f)$ and $H:=\Img(g)$. Then one gets the canonical exact sequence
$ 0\ra H\lraf{\epsilon_X} X\lraf{\lambda_X} {}_AU\ra 0$ of $A$-modules.
By Lemma \ref{module-A}(4), $\Img(m_X)=IX$. Let $i_X$ be the inclusion of $IX$ into $X$ and write $m_X=m'_Xi_X$. Consider the diagram
$$
\xymatrix{
&&
I\otimes_{A}U \ar@{->>}[d]^{m'_X} &&\\
&& IX \ar@{-->}[dl]_{m_{H}} \ar@{^(-_>}[d]^{i_X} &&\\
0\ar[r] & H \ar[r]^{\epsilon_{X}} & X \ar[r]^{\lambda_X} & U \ar[r] &0.}$$
As $\lambda_{X}$ is a homomorphism of $A$-modules, there holds $(IX)\lambda_{X}= I(\Img(\lambda_X))=IU=I \Coker(g)=0$ by Lemma \ref{module-A}(1). Thus there is an injective homomorphism $m_{H}: IX\to H$ of $A$-modules such that $i_X=m_H \epsilon_X$.  It follows from $_AIX\simeq {}_AI\otimes_A\Coker(g)$ that $m'_{X}$ is an isomorphism by Lemma \ref{noether}(2).

The isomorphism $N\otimes_B\Coker(f)\simeq \Img(g)/IX$  implies that $\theta_X$ is injective. Thus the proof of Lemma \ref{module-A}(3) implies that there is a  homomorphism $q: H\to N\otimes_BV$ making the following
diagram  commutative:
$$
\xymatrix@C=0.6cm@R=0.5cm{
& 0\ar[d] & 0\ar[d]\\
& IX\ar@{=}[r]\ar[d]^{m_{H}} & IX \ar[d]^-{i_X}\\
0 \ar[r] & H\ar[r]^{\epsilon_X}\ar[d]^-{q} & X\ar[r]^{\lambda_{X}}\ar[d]^-{p_X} & U\ar@{=}[d]\ar[r] & 0\\
0\ar[r] & N\otimes_{B}V\ar[r]^{\theta_{X}}\ar[d] & X/IX \ar[r]^-{\lambda'_{X}}\ar[d] & U\ar[r] & 0\\
& 0& 0
}
$$
This shows $\Coker(m_{H})\simeq N\otimes_{B}V$ as $A$-modules.
Now, we write $g=\sigma\epsilon_X$ with $\sigma: N\otimes_BY\ra H$ the canonical projection and $\epsilon_X: H\hookrightarrow X$ the inclusion, and consider the diagram
\begin{equation}\label{proof3.2}
\xymatrix{
N\otimes_{B}M\otimes_{A}X \ar@{->>}[d]_{1_{N\otimes_{B}M}\otimes \lambda_{X}}\ar[rr]^-{\psi\otimes 1_X} & & I\otimes_{A}X \ar[dl]^{1_I\otimes\lambda_X} \ar[dd]^{\mlt_X} &\\
N\otimes_{B}M\otimes_{A}U \ar[d]_-{1_N\otimes \eta_Y} \ar@{->>}[r]^-{\psi\otimes1_U}& I\otimes_{A}U \ar[d]^{m'_{X}m_H} & &\qquad \\
N\otimes_{B}Y \ar@{->>}[r]^{\sigma} & H\; \ar@{^(->}[r]^-{\epsilon_X} &X  && &}\end{equation}
Since $(X,Y,f,g)$ is a $\Lambda_{\psi}$-module, the out-side square in $(\ref{proof3.2})$ is commutative. By the definition of $m_H$, we know $\mlt_X=(1_I\otimes \lambda_X)(m'_Xm_H)\epsilon_X$. From the natural homomorphism $\psi: {}_AN\otimes_{B}M_A \ra {}_AI_A$, one sees that the upper square is commutative. As $1_{N\otimes_{B}M}\otimes \lambda_{X}$ is a surjective map and $\epsilon_X$ is an injective map, the down-left corner in $(\ref{proof3.2})$ is commutative.  This means that there is the following exact commutative diagram
$$
\xymatrix{
&N\otimes_{B}M\otimes_{A}U \ar[r]^-{1_N\otimes \eta_Y} \ar[d]^{\psi\otimes 1_U} & N\otimes_{B} Y \ar[d]^{\sigma}\ar[r]^{1_N\otimes \mu_Y} & N\otimes_{B} V \ar[r] \ar@{-->}[d]^-{\exists \, \nu }& 0\\
0\ar[r] & I\otimes_{A}U \ar[r]^{m'_{X}m_H} & H \ar[r]^-{q} &  N\otimes_{B} V \ar[r] & 0\\
}
$$
where the top row is exact because of $N\otimes_{B}-$ acting on the exact sequence ${}_BM\otimes_{A}U \lraf{\eta_{Y}} {}_BY \lraf{\mu_Y} {}_BV\lra 0$. By assumption, $A$ is left noetherian and $N\otimes_{B} V$ is a finitely generated $A$-module, whence every surjective endomorphism of $_AN\otimes_{B} V$ is an automorphism. Thus $\nu$ is an automorphism. This implies that $H$ is the pushout of $1_N\otimes_B \eta_Y$ and $\psi\otimes_A 1_U$. Since $MI=IU=0$, we have $I\otimes_AU\simeq I\otimes_{\Lambda}U $ and $N\otimes_{B}M\otimes_{A}U\simeq N\otimes_{B}M\otimes_{\Lambda}U$ as $A$-modules. Thus $H$ is  the pushout of  $1_N\otimes_B \eta_Y$ and $\psi\otimes_{\Lambda} 1_U$, as desired. $\square$

\medskip
For ${}_{\Lambda}X\in \Lambda\Modcat$, we define a quadruple $\mt_{\Lambda}(X):=(A\otimes_{\Lambda}X, {}_BM\otimes_{\Lambda}X, \pi_X,\psi_X)$, where
$\pi_X: M\otimes_AA\otimes_{\Lambda}X\simeq M\otimes_{\Lambda}X$ is the canonical homomorphism of $B$-modules, and where $$\psi_X: N\otimes_B(M\otimes_{\Lambda}X)\lra A\otimes_{\Lambda}X,\; n\otimes(m\otimes x)\mapsto (n\otimes m)\psi \otimes x \, \mbox{ for } n\in N, m\in M, x\in X,$$
is a homomorphism of $A$-modules. Note that $\psi_X=\Psi_{A\otimes_{\Lambda}X}$ and $\psi_{\Lambda}=\psi$. One can check that $(A\otimes_{\Lambda}X, {}_BM\otimes_{\Lambda}X, \pi_X,\psi_X)$ is a $\Lambda_{\psi}$-module. Moreover, for $\alpha\in \Hom_{\Lambda}(X,X')$, the pair $(1_A\otimes \alpha, 1_M\otimes\alpha)$ is a homomorphism from the $\Lambda_{\psi}$-module $\mt_{\Lambda}(X)$ to the one $\mt_{\Lambda}(X')$. Thus we get a functor

$$\mt_{\Lambda}: \Lambda\Modcat\lra \Lambda_{\psi}\Modcat, \quad {}_{\Lambda}X\mapsto \mt_{\Lambda}(X):=(A\otimes_{\Lambda}X, M\otimes_{\Lambda}X, \pi_X,\psi_X).$$
Note that $\mt_A(A\otimes_{\Lambda}X)=(A\otimes_{\Lambda}X, M\otimes_AA\otimes_{\Lambda}X, 1_{M\otimes_AA\otimes_{\Lambda}X}, \psi\otimes 1_{A\otimes X})\simeq \mt_{\Lambda}(X)$ as $\Lambda_{\psi}$-modules via the morphism $(1_{A\otimes_{\Lambda}X},\pi_{X}).$
In particular, $\mt_{\Lambda}(\Lambda)\simeq \Lambda_{\psi}e_1$, where $e_1=(\begin{smallmatrix}1_A &0\\ 0 & 0\end{smallmatrix})\in \Lambda_{\psi}$ with $1_A$ the identity of $A$.

To stress the $\Lambda$-decomposition of $_AA\otimes_{\Lambda}X$, we write sometimes $_AA\otimes_{\Lambda}X$ as $X(I):=X\oplus I\otimes_{\Lambda}X$.
Clearly, the $A$-module structure on $X(I)$ is given by $(\lambda, i)(x,j\otimes x')=(\lambda x, i\otimes x +\lambda j\otimes x')$, $\lambda\in\Lambda, i,j\in I$ and $x,x'\in X$. Thus $\mt_{\Lambda}(X)=(X(I), M\otimes_{\Lambda}X, \pi_X,\psi_X)).$

\begin{Lem}
\label{A-Lambda-homo}
Let $X$ and $X'$ be $\Lambda$-modules. Then there are the following homomorphisms (or isomorphisms) of abelian groups, which are natural in $X$ and $X'$.

$(1)$ $\Hom_{\Lambda}(X,X')\simeq \Hom_{A}(X(I), X')$, $f\mapsto \left(\begin{smallmatrix}
   f \\  0\\ \end{smallmatrix}\right)$ for $f\in \Hom_{\Lambda}(X,X')$.

$(2)$ $\Hom_{\Lambda}(X,I\otimes_{\Lambda}X')\ra  \Hom_{A}(X,X'(I))$, $g\mapsto (0,g)$ for $g\in \Hom_{\Lambda}(X,I\otimes_{\Lambda}X')$.

$(3)$ $\Hom_{\Lambda}(X, X'\oplus I\otimes_{\Lambda}X') \simeq \Hom_{A}(X(I),X'(I))$,$(a,c)\mapsto \left(\begin{smallmatrix}
           a & c \\ 0 & 1_{I}\otimes a \\ \end{smallmatrix}\right)$ for $a\in \Hom_{\Lambda}(X, X'), c\in\Hom_{\Lambda}(X, I\otimes_{\Lambda}X')$.
\end{Lem}

{\bf Proof.} $(1)$ There are isomorphisms $\Hom_{\Lambda}(X,X')\simeq
\Hom_{\Lambda}(X,\Hom_{A}(_{A}A_{\Lambda},X'))\simeq
\Hom_{A}(A\otimes_{\Lambda}X,X')\simeq \Hom_{A}(X(I),X')$
where the first isomorphism is induced from the isomorphism $X'\ra \Hom_{A}(_{A}A_{\Lambda},X'), \; x'\mapsto \{1_A \mapsto x'\}, x'\in X', $ the second is the adjunction  and the third is given by $X(I)\ra A\otimes_{\Lambda}X, (x,i\otimes y) \; \mapsto (1_{\Lambda},0)\otimes x+(0,i)\otimes y, x, y\in X.$
We can check that the composite of the above isomorphisms sends $f$ to $\binom{f}{0}$ for $f\in \Hom_{\Lambda}(X,X')$, and is natural in each variables.

$(2)$ For $g\in \Hom_{\Lambda}(X,I\otimes_{\Lambda}X')$, it suffices to prove that
$(0,g): X\to X'(I)$ is a homomorphism of $A$-modules. In fact, take $a=(\lambda,i)\in A$, $x\in X$ and consider the $\Lambda$-module $X$ as $A$-module, that is $ix=0$ and $ax=(\lambda, i)x= \lambda x$. Then $(ax)(0,g)=(0,(ax)g)=(0,(\lambda x)g)$. On the other hand, $a(0,(x)g)=(\lambda,i)(0,(x)g)=(0,\lambda(x)g)=(0,(\lambda x)g)$. Thus $(0,g)$ is a homomorphism of $A$-modules.

$(3)$ The proof is similar to the one of (1). Since we have the isomorphisms  $\Hom_{\Lambda}(X, X'\oplus I\otimes_{\Lambda}X') \simeq \Hom_{\Lambda}\big(X,\Hom_A(_AA_{\Lambda},X'(I))\big)\simeq
\Hom_{A}\big(A\otimes_{\Lambda}X,X'(I)\big)
\simeq \Hom_{A}\big(X(I),X'(I)\big)$, their composite sends $(a,c)$ to
$\left(\begin{smallmatrix} a & c \\ 0 & 1_{I}\otimes a \\ \end{smallmatrix}\right)$ for $a\in
\Hom_{\Lambda}(X, X')$, $c\in\Hom_{\Lambda}(X,
           I\otimes_{\Lambda}X')$. All isomorphisms are natural in each variables. $\square$

\smallskip
The next lemma for $\Lambda_{(0,0)}$ was indicated in \cite{GaoP17}. We  will state it for $\Lambda_{\psi}$ and include more details for our applications. Note that the functor $\mz_B: B\Modcat\to \Lambda_{\psi}\Modcat$, $_BY\mapsto (0,Y,0,0)$ is well defined.

\begin{Lem} \label{homo-equ}
Let $_{\Lambda}X, X'\in \Lambda\modcat$ and $_BY$, $_BY'\in B\modcat$. Then there are the following isomorphisms which are natural in each variables.	

$(1)$ $\Hom_{\Lambda}(X, N\otimes_BY) \simeq \Hom_{\Lambda_{\psi}}\big(\mt_{\Lambda}(X), \mt_B(Y)\big)$, $f\mapsto (\left(\begin{smallmatrix}   f \\    0\\ \end{smallmatrix}\right),0)$.
			
$(2)$ $\Hom_{\Lambda}(X, X') \simeq  \Hom_{\Lambda_{\psi}}\big(\mt_{\Lambda}(X), \mz_{\Lambda}(X')\big) $, $g\mapsto (\left(\begin{smallmatrix}    g \\    0\\ \end{smallmatrix}\right),0)$.

$(3)$ $\Hom_{\Lambda}(X, X'\oplus I\otimes_{\Lambda}X') \simeq  \Hom_{\Lambda_{\psi}}\big(\mt_{\Lambda}(X), \mt_{\Lambda}(X')\big) $, $(a,c)\mapsto \big(\left(\begin{smallmatrix}   a & c \\   0 & 1_{I}\otimes a \\ \end{smallmatrix}\right),1_{M}\otimes a \big)$ for $a: X\to X'$ and $c: X\to I\otimes_{\Lambda}X'$.

$(4)$ $ \Hom_B(Y, M\otimes_{\Lambda}X) \simeq \Hom_{\Lambda_{\psi}}\big(\mt_{B}(Y), \mt_{\Lambda}(X)\big)$, $h\mapsto ((1_N\otimes h)\psi_{X},h)$.
			
$(5)$  $\Hom_B(Y, Y') \simeq \Hom_{\Lambda_{\psi}}\big(\mt_{B}(Y), \mt_B(Y')\big)$, $t \mapsto (1_N\otimes t,t) $.			

$(6)$ $ \Hom_B(Y, Y')    \simeq \Hom_{\Lambda_{\psi}}\big(\mt_{B}(Y), \mz_{B}(Y')\big) $, $t \mapsto (0, t)$.
		
$(7)$  $ \Hom_{\Lambda_{\psi}}(\mt_{B}(Y), \mz_{\Lambda}(X))= \Hom_{\Lambda_{\psi}}\big(\mt_{\Lambda}(X), \mz_{B}(Y)\big)=0$.
\end{Lem}

{\bf Proof.}
$(1)$  There are isomorphisms $\Hom_{\Lambda}(X, N\otimes_BY) \simeq  \Hom_{A}(X(I), N\otimes_BY) \simeq  \Hom_{\Lambda_{\psi}}(\mt_{A}(X(I)), \mt_B(Y))$ $\simeq \Hom_{\Lambda_{\psi}}(\mt_{\Lambda}(X), \mt_B(Y))$
where the first isomorphism is given by Lemma \ref{A-Lambda-homo}(3), the second one follows from the adjoint pair $(\mt_{A},\mU_A)$ of functors in  Lemma \ref{mr-recol} and $\mU_{A}\mt_{B}(Y)= N\otimes_{B}Y$, and the third is given by the isomorphism $\mt_{A}(X(I))\simeq \mt_{\Lambda}(X)$.  Verifications show that the composite of theses isomorphisms sends $f\in \Hom_{\Lambda}(X, N\otimes_BY)$ to $(\binom{f}{0},0)$ and is natural in each variables.

$(2)$ Its proof is similar to the one of (1).

$(3)$ Due to Lemma \ref{A-Lambda-homo}(3), we have $\Hom_{\Lambda}(X,  X'\oplus I\otimes_{\Lambda}X') =  \Hom_{A}(X(I), X'(I)).$ Since $\mt_{A}$ is fully faithful by Lemma \ref{A-Lambda-homo}(1), we have $\Hom_{A}(X(I), X'(I))\simeq  \Hom_{\Lambda_{\psi}}(\mt_{A}(X(I)), \mt_{A}(X'(I))).$  Due to $\mt_{A}(X(I))\simeq \mt_{\Lambda}(X)$ for any $X$, we get $\Hom_{\Lambda_{\psi}}(\mt_{A}(X(I)), \mt_{A}(X'(I)))~\simeq \Hom_{\Lambda_{\psi}}(\mt_{\Lambda}(X), \mt_{\Lambda}(X'))$. By verification, the composite of these isomorphisms sends $(a,c)\in \Hom_{\Lambda}(X, X'\oplus I\otimes_{\Lambda}X')$ to $(\left(\begin{smallmatrix}   a & c \\   0 & 1_{I}\otimes a \\ \end{smallmatrix}\right),1_{M}\otimes a )$, and is natural in each variable.

$(4)$ It follows from the adjoint pair $(\mt_{B},\mU_B)$ in Lemma \ref{mr-recol} and $\mU_{B}\mt_{\Lambda}(X) = M\otimes_{\Lambda}X$ that $\Hom_{B}(Y, M\otimes_{\Lambda}X) \simeq \Hom_{\Lambda_{\psi}}(\mt_B(Y),\mt_{\Lambda}(X)),$ $h\mapsto ((1_N\otimes h)\psi_{X},h)$, which is natural in each variable.

$(5)$ This is a consequence of the fully faithful functor $\mt_{B}$ in Lemma \ref{mr-recol}.

(6) The proof is similar to the one of (4).

$(7)$  By $\mU_{B}\mz_{\Lambda}(X)=0$ and the adjoint pair $(\mt_B,\mU_B)$ of functors, we have $\Hom_{\Lambda_{\psi}}(\mt_{B}(Y), \mz_{\Lambda}(X)) \simeq \Hom_B(Y, \mU_B\mz_{\Lambda}(X))$ = $0$.
Similarly, since $(\mt_{A},\mU_A)$ is an adjoint pair of functors and $\mU_{A}\mz_{B}(Y)=0$, we have $\Hom_{\Lambda_{\psi}}(\mt_{\Lambda}(X), $ $\mz_{B}(Y))$ $\simeq \Hom_{\Lambda_{\psi}}(\mt_{A}(X(I)), \mz_{B}(Y))\simeq \Hom_A(X(I), \mU_A\mz_B(Y))=0$. $\square$

\begin{Theo} \label{main-result-1}Suppose that
$_\Lambda N_B$, $_BM_\Lambda$, $_{\Lambda}I_{\Lambda}$ are weakly compatible bimodules and $(X,Y,f,g)$ is a $\Lambda_{\psi}$-module. Then $(X,Y,f,g)$ is Gorenstein-projective if the following hold:

$(a)$ Both $_{\Lambda}\Coker(g)$ and $_B\Coker(f)$ are Gorenstein-projective, and	

$(b)$	${}_BM\otimes_A\Coker(g)\simeq {}_B\Img(f)$, $_AN\otimes_B\Coker(f)\simeq {}_A\Img(g)/IX$, and
${}_AI\otimes_A\Coker(g)\simeq {}_AIX$, where $\Coker(f)$ and $\Img(g)$ denote the cokernel of $f$ and the image of $g$, respectively.
\end{Theo}

{\bf Proof.} Suppose that $(a)$ and $(b)$ hold true.
Then ${}_{\Lambda}U:=\Coker(g)$ and ${}_BV:=\Coker(f)$ are Gorenstein-projective. By definition, there are two totally exact sequences of projective modules
$$ P^{\bullet}:\quad 	\cdots \lra P^{-1}\lraf{d_P^{-1}}P^0\lraf{d_P^0}P^1\lra\cdots, \mbox{ and } \; Q^{\bullet}:\quad	\cdots \lra Q^{-1}\lraf{d_Q^{-1}}Q^0\lraf{d_Q^0}Q^1\lra\cdots $$
over $\Lambda$ and $B$, respectively, such that $_{\Lambda}U=\Ker(d_{P}^0)$ and $_BV=\Ker(d_{Q}^0)$. To prove $(X,Y,f,g)$ is a Gorenstein-projective, we shall construct a totally exact complex $\cpx{T}=(T^i,d_T^i)_{i\in \mathbb{Z}}\in \mathscr{C}(\pmodcat{\Lambda_{\psi}} )$, such that $\Ker(d_T^0)\simeq (X,Y,f,g)$ as $\Lambda_{\psi}$-modules. We define $T^i:=\mt_{\Lambda}(P^i)\oplus \mt_B(Q^i)$ for all $i\in \mathbb{Z}$. Then $T^i\in \mathscr{C}(\pmodcat{\Lambda_{\psi}})$ by Lemma \ref{prop:projmod}(1). To define $ d^i_T$, we have to define a few families of maps.

(1) A homomorphism $\rho^i: Q^i\to M\otimes_{\Lambda}P^{i+1}$ of $B$-modules for $i\in\mathbb{Z}$.

By assumptions, $_BM_{\Lambda}$ is weakly compatible, this implies that $M\otimes_{\Lambda}\cpx{P}$ is exact by Definition \ref{weak}(C3), and therefore $\Ker(1_M\otimes d^0_P)=M\otimes_{\Lambda}U$.  Note that $M\otimes_{\Lambda}P^i\in \add(_BM)$ for all $i$. Thus it follows from Definition \ref{weak}(C1) that $\Ext^{1}_{B}(\Ker(d_{Q}^{i}), M\otimes_{\Lambda}P^{i})=0$ for all $i\ge 0$ and $\Ext^1_B(Q^{-i}, \Img(1_M\otimes d_P^{-i}))=0$ for $i\ge 1$. Thus, starting with the exact sequence $$ (*)\quad 0\lra {}_BM\otimes_{\Lambda}U \lraf{\eta_{Y}} {}_BY \lraf{\mu_Y} {}_BV\lra 0$$ and applying Horseshoe Lemma \ref{horseshoe}, we get an exact sequence of complex $$(**)\quad 0\lra M\otimes_{\Lambda} \cpx{P} \lraf{\cpx{a}} \cpx{Y} \lraf{\cpx{b}} \cpx{Q}\lra 0$$
of $B$-modules, such that $Y^{i}:=M\otimes_\Lambda P^{i}\oplus Q^{i} $ and $d_{Y}^{i}=
(\begin{smallmatrix}
	 	1_M\otimes d_P^{i} & 0 \\
	 	\rho^{i} & d_ Q^{i}
	 \end{smallmatrix})$, where $\rho^{i}: Q^{i}\ra M\otimes_{\Lambda}P^{i+1}$ is a homomorphism of $B$-modules, and where
$a^i=(1,0)$ and $b^{i}=(\begin{smallmatrix}
	 	 0 \\
	 	1
	 \end{smallmatrix})$ are canonical maps for all $i$. Note that by taking kernels of $(**)$ at degree $0$, we get back the exact sequence $(*)$.

(2) Two homomorphisms $\alpha^{i}: P^{i}\ra I\otimes_{\Lambda}P^{i+1}$ and $\beta^{i}: P^{i}\ra N\otimes_{B}Q^{i+1}$ of $\Lambda$-modules for all $i\in \mathbb{Z}$.

In fact, we write $H:=\Img(g)$ and get a canonical exact sequence
$ 0\ra H\ra X\ra {}_AU\ra 0$ of $A$-modules,
which restricts to an exact sequence of $\Lambda$-modules
$$ (\dag)\quad 0\lra {}_{\Lambda}H\lraf{\epsilon_X} {}_{\Lambda}X\lraf{\lambda_{X}}
{}_{\Lambda}U\lra 0.$$
Now, we define
$Z^i:=(I\otimes_\Lambda P^{i})\oplus (N\otimes_BQ^{i}) $, $d_{Z}^{i}=\bigl(\begin{smallmatrix}
 1_I \otimes d^i_P & 0 \\
 \tau^{i} & 1_N \otimes d_Q^{i}
\end{smallmatrix}\bigr)$, where $\tau^{i}=(1_N \otimes \rho^i)(\psi \otimes 1_{P^{i+1}})$ is the composite of the homomorphisms of $\Lambda$-modules:
$$ {}_{\Lambda}N\otimes_{B}Q^{i} \lraf{1_N \otimes \rho^i} {}_{\Lambda}N\otimes_BM\otimes_{\Lambda} P^{i+1}\lraf{\psi \otimes 1_{P^{i+1}}} {}_{\Lambda}I \otimes_{\Lambda} P^{i+1}.$$ We shall show that $\cpx{Z}=(Z^i,d_Z^i)_{i\in\mathbb{Z}}$ is an exact complex such that $H\simeq \Ker(d^0_Z)$. Indeed, it follows from the complex $\cpx{Y}$ that $$\begin{pmatrix}
	 	1_M\otimes d_P^{i} & 0 \\
	 	\rho^{i} & d_ Q^{i}
	 \end{pmatrix}\begin{pmatrix}
	 	1_M\otimes d_P^{i+1} & 0 \\
	 	\rho^{i+1} & d_ Q^{i+1}
	 \end{pmatrix}=\begin{pmatrix}
	 	0 & 0 \\
	 	\rho^{i}(1_M\otimes d_P^{i+1})+ d_ Q^{i}\rho^{i+1} & 0
	 \end{pmatrix}=0.$$
This implies that $\rho^{i}(1_M\otimes d_P^{i+1})+ d_ Q^{i}\rho^{i+1}=0$ and
$ (1_N\otimes \rho^{i})(1_N\otimes1_M\otimes d_P^{i+1})+ (1_N\otimes d_ Q^{i})(1_N\otimes\rho^{i+1})=0.$
By multiplying $\psi\otimes 1_{P^{i+2}}$, we further obtain
$$ (1_N\otimes \rho^{i})(1_N\otimes1_M\otimes d_P^{i+1})(\psi\otimes 1_{P^{i+2}})+ (1_N\otimes d_ Q^{i})(1_N\otimes\rho^{i+1})(\psi\otimes 1_{P^{i+2}})=0.                              $$
Due to $(1_N\otimes1_M\otimes d_P^{i+1})(\psi\otimes 1_{P^{i+2}})=(\psi \otimes 1_{P^{i+1}})(1_I\otimes d_P^{i+1})$, we get
$ \tau^{i}(1_I\otimes d_P^{i+1})+ (1_N\otimes d_ Q^{i})\tau^{i+1}=0.$ Thus
$$ d^i_Zd^{i+1}_Z=\begin{pmatrix}
 1_I \otimes d^i_P & 0 \\
 \tau^{i} & 1_N \otimes d_Q^{i}
\end{pmatrix} \begin{pmatrix}
 1_I \otimes d^{i+1}_P & 0 \\
 \tau^{i+1} & 1_N \otimes d_Q^{i+1}
\end{pmatrix}=0,$$
and $\cpx{Z}$ is a complex. Moreover, there is the exact sequence of complexes of $\Lambda$-modules
$$0\lra I\otimes_{\Lambda} \cpx{P} \lraf{\cpx{c}} \cpx{Z} \lraf{\cpx{d}} N\otimes_B \cpx{Q}\lra 0,$$
where $\cpx{c}$ and $\cpx{d}$ are canonical inclusion and projection, respectively. By the weak compatibility of $_{\Lambda}I_{\Lambda}$ and $_{\Lambda}N_B$, both $ I\otimes_{\Lambda} \cpx{P}$ and $N\otimes_B \cpx{Q}$ are exact complexes of $\Lambda$-modules.
Thus $\cpx{Z}$ is also an exact complex.

We now show $\Ker(d^0_Z)\simeq H.$
Let $\sigma^i:=\bigl(\begin{smallmatrix}
  \psi \otimes 1_{P^{i}} & 0 \\
 0 & 1_{N}\otimes 1_{Q^i}
\end{smallmatrix}\bigr)_{i\in Z}$ and $\cpx{\sigma}=(\sigma^i)_{i\in \mathbb{Z}}$. Then $\cpx{\sigma}: N\otimes_B\cpx{Y}\ra \cpx{Z}$ is a chain map of complexes such that the following is an exact commutative diagram in $\mathscr{C}(A)$:
$$
\xymatrix{
&N\otimes_{B}M\otimes_{\Lambda}\cpx{P} \ar[r]^-{1_N\otimes \cpx{a}} \ar[d]^{\psi\otimes 1_\cpx{P}} & N\otimes_{B}\cpx{Y} \ar[d]^{\cpx{\sigma}}\ar[r]^{1_N\otimes \cpx{b}} & N\otimes_{B}\cpx{Q} \ar[r] \ar@{=}[d]& 0\\
0\ar[r] & I\otimes_{\Lambda}\cpx{P}\ar[r]^{\cpx{c}} & \cpx{Z} \ar[r]^{\cpx{d}} &  N\otimes_{B}\cpx{Q} \ar[r] & 0.  }
$$
Consider the differentials in degree $0$, there is the commutative diagram with exact rows:

$$\xymatrix@C=0.6cm@R=0.6cm{
  & N\otimes M\otimes P^1 \ar[rr]^-{1_N\otimes a^1} \ar'[d]^(0.53){\psi\otimes 1_{P^1}}[dd]
      &  & N\otimes Y^1 \ar'[d]^{\sigma^1}[dd]  \ar[rr]^-{1_N\otimes b^1} & & N\otimes_BQ^1 \ar@{=}[dd]     \\
  N\otimes M\otimes P^0 \ar[ur]^{1_N\otimes 1_M\otimes d_P^0}\ar[rr]_(0.7){1_N\otimes a^0}\ar[dd]^(0.33){\psi\otimes 1_{P^0}}
      &  & N\otimes Y^0 \ar[ur]^{1_N\otimes d_Y^0}\ar[dd]^(0.33){\sigma^0} \ar[rr]_(0.7){1_N\otimes b^0} & &N\otimes_B Q^0 \ar[ur]^{1_N\otimes d_Q^0}\ar@{=}[dd] \\
  & I\otimes P^1 \ar'[r]^{c^1}[rr]
      &  & Z^1    \ar'[r]^{d^1}[rr] & &    N\otimes Q^1         \\
  I\otimes P^0 \ar[rr]^-{c^0}\ar[ur]^{1_I\otimes d^0_P}
      &  & Z^0 \ar[ur]^{d_Z^0}     \ar[rr]^{d^0} & & N\otimes_B Q^0 \ar[ur]^{1_N\otimes d_Q^0}  }
$$
By enlarging the diagram forward, we get the exact commutative diagram of $A$-modules by Lemma \ref{ker-seq-sur}(2)
$$
\xymatrix{
&N\otimes_{B}M\otimes_{\Lambda}U \ar[r]^-{1_N\otimes \eta_Y} \ar[d]^{\psi\otimes 1_U} & N\otimes_{B} Y \ar[d]^{\delta}\ar[r]^{1_N\otimes \mu_Y} & N\otimes_{B} V \ar[r] \ar@{=}[d]& 0\\
0\ar[r] & I\otimes_{\Lambda}U \ar[r]^{c^{0}} & \Ker(d_{Z}^0) \ar[r]^{d^0} &  N\otimes_{B} V \ar[r] & 0 }
$$
This shows that $\Ker(d_Z^0)$ is the pullback of $\psi\otimes_{\Lambda} 1_U$ and $1_N\otimes_B\eta_Y$.
By Lemma \ref{b-condition}, $H$ is also the pushout of $\psi\otimes_{\Lambda} 1_U$ and $1_N\otimes_{B}\eta_Y$. Thus $H\simeq \Ker(d_Z^0)$ as $A$-modules.
Thus there is the decomposition of $d_Z^{-1}$
$$\xymatrix{\cpx{Z}: \quad \cdots \ar[r]^-{d_Z^{-2}}  &(I\otimes_\Lambda P^{-1})\oplus (N\otimes_BQ^{-1})\ar@{->>}[dr] \ar[rr]^{d^{-1}_Z} && (I\otimes_\Lambda P^{0})\oplus (N\otimes_BQ^{0})\ar[r]^-{d_Z^0} & \cdots \\
& & H \ar[ur]^-{d_H} & &}
$$
so that the following two diagrams are commutative
\begin{equation}
\xymatrix{
0 \ar[r] & I\otimes_{\Lambda}U \ar[r]^-{m'_{X}m_{H}} \ar[d]^{1_I\otimes d_{U}} & H \ar[d]^{d_{H}}\ar[r]^-{q} & N\otimes_{B} V \ar[r] \ar[d]& 0  \\
0\ar[r] & I\otimes_{\Lambda}P^{0} \ar[r]^-{(1,0)} &  I\otimes_{\Lambda}P^{0} \oplus N\otimes_{B} Q^{0} \ar[r]^-{\binom{0}{1}} &  N\otimes_{B} Q^{0} \ar[r] & 0, & }
\end{equation}
$$
\xymatrix{
 & N\otimes_{B}Y \ar[r]^-{1_N \otimes d_Y} \ar[d]^{\sigma} &   N\otimes_{B}Y^0 \ar[d]^{\sigma^0}\ar[r]^-{1_N \otimes d_Y^0} & N\otimes_{B} Y^1\ar[d]\ar[r]\ar[d]^-{\sigma^1} &\cdots \\
0\ar[r] & H \ar[r]^-{d_H} &  Z^{0} \ar[r]^-{d_Z^0} &  Z^{1} \ar[r] &\cdots & }
$$

Finally, we give the definition of $\alpha^i$ and $\beta^i$. In fact, there is the exact sequence $$(\dag)\quad 0\lra {}_{\Lambda}H\lraf{\epsilon_X} {}_{\Lambda}X\lraf{\lambda_{X}}
{}_{\Lambda}U\lra 0.$$
Since $_{\Lambda}I_{\Lambda}$ is weakly compatible, $\Hom_{\Lambda}(\cpx{P}, {}_{\Lambda}I)$ is exact. This yields $\Ext^1_{\Lambda}(\Ker(d^i_P), I)=0$ for all $i$. Due to $I\otimes_{\Lambda}P^i\in \add(_{\Lambda}I)$, there holds $\Ext^1_{\Lambda}(\Ker(d^i_P), I\otimes_{}P^i)=0$ for all $i\ge 0$. Similarly, the weak compatibility of $_{\Lambda}N_B$ implies that $\Ext^1_{\Lambda}(\Ker(d^i_P),N\otimes_BQ^i)=0$ for all $i$. Thus $\Ext^1_{\Lambda}(\Ker(d_P^i),Z^i)=0$ for all $i\ge 0$. Clearly, $\Ext^1_{\Lambda}(P^{-i}, \Img(d_Z^{-i})=0$ for all $i\ge 1$. Hence, applying Lemma \ref{horseshoe} to the exact sequence $(\dag)$ and the exact complexes $\cpx{Z}$ and $\cpx{P}$, we get an exact sequence of exact complexes
$$(\dag\dag)\quad  0\lra \cpx{Z}\lraf{\cpx{p}}\cpx{E}\lraf{\cpx{q}}\cpx{P}\lra 0 $$
in $\mathscr{C}(\Lambda),$ where $\cpx{p}$ and $\cpx{q}$ are canonical inclusion and projection, respectively, and where $E^i:= P^i \oplus (I\otimes_{\Lambda}P^i)\oplus(N\otimes_BQ^{i})\in \Lambda \modcat$, and
$$d_E^i=
\begin{pmatrix}
 d_{P}^{i} & \alpha^{i} & \beta^{i} \\
 0 & 1_{I}\otimes d^i_P & 0\\
 0 & \tau^{i}& 1_{N}\otimes d^i_Q
\end{pmatrix}: P^i \oplus (I\otimes_{\Lambda}P^i)\oplus(N\otimes_BQ^{i})\lra P^{i+1} \oplus (I\otimes_{\Lambda}P^{i+1})\oplus(N\otimes_BQ^{i+1}),$$ $\alpha^{i}: P^{i}\ra I\otimes_{\Lambda}P^{i+1}$, and $\beta^{i}: P^{i}\ra N\otimes_{B}Q^{i+1}$ are all homomorphisms of $\Lambda$-modules. Recall that $\tau^{i}=(1_N \otimes \rho^i)(\psi \otimes 1_{P^{i+1}}): N\otimes_BQ^i\ra I\otimes_{\Lambda}P^{i+1}$ is a homomorphism of $\Lambda$-modules. Here, compared with Lemma \ref{horseshoe}, the order of direct summands of $E^i$ is changed. Note that $(\dag)$ is gotten by taking kernels at degree $0$ in $(\dag\dag)$. Visually, the positive part of $(\dag\dag)$ looks as follows:
\begin{equation}\label{diag6}
\xymatrix{
   & 0 \ar[d] & 0 \ar[d] & 0 \ar[d] & \\
  0 \ar[r] &  H \ar[d]_{\epsilon_X} \ar[r]^-{d_H} & I\otimes_{\Lambda}P^{0}\oplus N\otimes_{B}Q^{0} \ar[d]^{\bigl(\begin{smallmatrix}
 0 & 1 & 0\\ 0 & 0& 1 \\
\end{smallmatrix}\bigr)}
  \ar[r]^{d_Z^0} & I\otimes_{\Lambda}P^{1}\oplus N\otimes_{B}Q^{1} \ar[d]^{\bigl(\begin{smallmatrix}
 0 & 1 & 0\\ 0 & 0& 1 \\
\end{smallmatrix}\bigr)} \ar[r] & \cdots  \\
 0 \ar@{-->}[r] &  {}_{\Lambda}X \ar@{-->}[ur]^-{(e_1,e_2)} \ar@{-->}[r]_-{d_{_{\Lambda}X}}  \ar[d]_{\lambda_{X}} & P^0 \oplus I\otimes_{\Lambda}P^{0}\oplus N\otimes_{B}Q^{0} \ar@{-->}[r]_{d_E^0} \ar[d]^{\bigl(\begin{smallmatrix}
 1 \\ 0\\ 0\\
\end{smallmatrix}\bigr)} &  P^1 \oplus I\otimes_{\Lambda}P^{1}\oplus N\otimes_{B}Q^{1}  \ar@{-->}[r] \ar[d]^{\bigl(\begin{smallmatrix}
 1 \\ 0\\ 0\\
\end{smallmatrix}\bigr)} & \cdots &  \\
 0 \ar[r] &  U \ar[d] \ar[r]^-{d_U} & P^0 \ar[d] \ar[r]^{d^0_P} & P^1 \ar[d] \ar[r] & \cdots  \\
  & 0 & 0 & 0 &  }
\end{equation}
where $d_{_{\Lambda}X}=(\lambda_{X}d_{U},e_1,e_2)$, both $e_1: {}_{\Lambda}X \ra I\otimes_{\Lambda}P^{0}$ and $e_2: {}_{\Lambda}X \ra N\otimes_{B}Q^{0}$ are homomorphisms of $\Lambda$-modules.

Observe that $1_N\otimes d^i_Q$ and $(0,\tau^i)$ are homomorphisms of $A$-modules. The term $P^i\oplus I\otimes_{\Lambda}P^i$ has an $A$-module structure which is isomorphic to $_AA\otimes_{\Lambda}P^i$. By Lemma \ref{A-Lambda-homo}(3), the map
$$\begin{pmatrix}
 d_{P}^{i} & \alpha^{i}  \\
 0 & 1_{I}\otimes d^i_P\\
\end{pmatrix}: P^i \oplus I\otimes_{\Lambda}P^i \lra P^{i+1} \oplus I\otimes_{\Lambda}P^{i+1} $$
is the image of the homomorphism $(d_P^i, \alpha^i): P^i\to P^{i+1}\oplus I\otimes_{\Lambda}P^{i+1}$ of $\Lambda$-modules, thus it is a homomorphism of $A$-modules. Similarly, $\big(\begin{smallmatrix}\beta^i\\ 0 \end{smallmatrix}\big): P^i(I)\ra N\otimes_BQ^{i+1}$ is an $A$-module homomorphism.
Let $$F^i:=P^i(I)\oplus N\otimes_BQ^i, \quad d_F^i:=\begin{pmatrix}
 d_P^i & \alpha^i & \beta^i \\
 0 & 1_I\otimes d^i_P &0\\
 0 & \tau^i &1_N\otimes d^i_Q\\
\end{pmatrix}.$$ Then $\cpx{F}=(F^i, d^i_F)$ is a complex of $A$-modules, which is exact since the restriction of $\cpx{F}$ to $\Lambda$-modules is the exact complex $\cpx{E}$.

We show $\Ker(d_F^0)\simeq X$ as $A$-modules. First, the exact sequence $ 0\ra {}_{\Lambda}X \lraf{d_{_{\Lambda}X}} E^{0} \lraf{d_{E}^{0}} E^{1} $ of $\Lambda$-modules, with $d_{_{\Lambda}X}=(\lambda_{X}d_{U},e_1,e_2): X\ra P^0 \oplus I\otimes_{\Lambda}P^0\oplus N\otimes_BQ^0$, can be regarded as an exact sequence of $\mathbb{Z}$-modules:
$$ 0\lra X \lraf{d_X} F^{0} \lraf{d_{F}^{0}} F^{1} $$
where $d_{X}:=\big((\lambda_{X}d_{U},e_1),e_2\big):X\ra F^{0}=(P^{0}\oplus (I\otimes_{\Lambda}P^{0}))\oplus (N\otimes_BQ^{1})$.
It suffices to show that $d_{X}$ is a homomorphism of $A$-modules.

Let $x\in X$, $a=(\lambda,i)\in A=\Lambda\ltimes I$. Then $(ax)d_{X}=\big(((ax)\lambda_{X}d_U,(ax)e_1),(ax)e_2\big)$.
Since $\lambda_{X}$ is an $A$-homomorphism, it follows from Lemma \ref{module-A} that $(IX)\lambda_{X}=0$ and
$ (ax)\lambda_{X}d_U=(\lambda x)\lambda_{X}d_U=\lambda(x)\lambda_{X}d_U.$ Note that $(ix)(e_1,e_2)=(i\otimes x)\mlt_{X}(e_1,e_2)$. We deduce from the diagrams $(6), (4)$ and $(5)$ that
$\mlt_{X}(e_1,e_2) = (1_I \otimes \lambda_{X})m'_{X}m_{H}\epsilon_{X}(e_1,e_2)$ = $(1_I \otimes \lambda_{X})m'_{X}m_{H}d_H = (1_I \otimes \lambda_{X})(1_I \otimes d_U)(1,0) = (1_I \otimes \lambda_{X}d_U,0)$. Thus $(ix)(e_1,e_2)=\big(i\otimes (x)\lambda_{X}d_U,0\big).$
As $e_i$ is a $\Lambda$-homomorphism,  there holds $(\lambda x)e_i = \lambda(x)e_i$, $i=1,2$. Hence
$$ (ax)d_X=\big((\lambda(x)\lambda_{X}d_U,   \lambda(x)e_1 +i\otimes (x)\lambda_{X}d_U),\lambda(x)e_2\big) .$$
It then follows from the module structure of $P^{0}(I)$ and $N\otimes_{B} Q^{0}$ that
$ (\lambda,i)\big((x)\lambda_{X}d_U,(x)e_1\big) = \big(\lambda(x)\lambda_{X}d_U, \lambda(x)e_1 + i\otimes (x)\lambda_{X}d_U\big)$
and $(\lambda,i)(x)e_2=\lambda (x)e_2.$ This shows
$$\begin{array}{rl} a(x)d_X & =(\lambda,i)(x)d_{X}=(\lambda,i)\big(((x)\lambda_{X}d_U,(x)e_1),(x)e_2\big)\\
 & = \big((\lambda(x)\lambda_{X}d_U,   \lambda(x)e_1 +i\otimes (x)\lambda_{X}d_U),\lambda(x)e_2\big)=(ax)d_{X},\end{array}$$
that is, $d_X$ is an $A$-homomorphism. Thus $_AX$ has a complete projective resolution $\cpx{F}$.

\medskip
At this stage, we can define the complex $\cpx{T}=(T^i,d_T^i)$ of $\Lambda_{\psi}$-modules as follows. Let
$$T^i:=\mt_{\Lambda}(P^i)\oplus \mt_{B}(Q^i)= \bigl(F^i, Y^i, \big(\begin{smallmatrix}\pi_{P^i} & 0\\ 0 &0 \end{smallmatrix}\big), \big(\begin{smallmatrix}\psi_{P^i} & 0\\ 0 &1_{N\otimes_BQ^i} \end{smallmatrix}\big)\bigr),
\; d^i_T:=(d_F^i,d_Y^i)=\big(\big(\begin{smallmatrix}
 d_{P}^{i} & \alpha^{i} & \beta^i \\
 0 & 1_{I}\otimes d^i_P & 0\\
 0 &\tau^i & 1_N\otimes d_Q^i\\
\end{smallmatrix}\big), \big(\begin{smallmatrix}1_M\otimes d_P^i & 0\\ \rho^i & d^i_Q \end{smallmatrix}\big)\big).$$
Note that $\Phi_{Q^i}=0$ since we assume $\phi=0$ in the Morita context ring $\Lambda_{\psi}.$

We have to check that $(d_F^i,d_Y^i)$ is a homomorphism of $\Lambda_{\psi}$-modules. In fact, $(d_F^i,d_Y^i)$ can be written as $ \begin{pmatrix}
 t_{11}^{i} & t_{12}^{i}  \\
 t_{21}^{i} & t_{22}^{i}
\end{pmatrix}$
with $$t_{11}^i = \big((\begin{smallmatrix}
 d_{P}^{i} & \alpha^{i}  \\
 0 & 1_{I}\otimes d^i_P
\end{smallmatrix}),1_{M}\otimes d^i_P\big): \mt_{\Lambda}(P^{i}) \lra \mt_{\Lambda}(P^{i+1}), \qquad t_{12}^{i}=\big(\left(\begin{smallmatrix}
   \beta^{i} \\
   0
\end{smallmatrix}\right),0): \mt_{\Lambda}(P^{i})\lra \mt_B(Q^{i+1}),$$
$$t_{21}^{i}=((0,\tau^{i}),\rho^{i}):  \mt_{B}(Q^{i})\lra \mt_{\Lambda}(P^{i+1}), \qquad t_{22}^{i}=(1_N\otimes d_{Q}^{i},d_{Q}^{i}): \mt_{B}(Q^{i})\lra \mt_B(Q^{i+1}).$$
To see that these $t_{pq}^i$ are homomorphisms of $\Lambda_{\psi}$-modules, we just note that $t_{11}^{i}$ is the image of $(d_{P}^{i}, \alpha^{i})$ under the isomorphism in Lemma \ref{homo-equ}(3). Similarly, it follows from Lemmas \ref{homo-equ}(1),(4) and (5) that $t^i_{12}, t^i_{21}$ and $t^i_{22}$ are homomorphisms of $\Lambda_{\psi}$-modules.

We show that $\cpx{T}$ is a total projective resolution of $(X,Y,f,g)$.

Actually, the complex $\cpx{T}$ is exact because $\cpx{F}$ and $\cpx{Y}$ are exact. By Lemma \ref{prop:projmod}, each term of $\cpx{T}$ is a projective $\Lambda_{\psi}$-module. Thus $\cpx{T}$ is an exact complex in $\mathscr{C}(\pmodcat{\Lambda_{\psi}})$.

Next, we show that the complex $\Hom_{\Lambda_{\psi}}(\cpx{T},\Lambda_{\psi})$ is exact. This is equivalent to saying that the complexes $\Hom_{\Lambda_{\psi}}(\cpx{T}, \mt_{\Lambda}(\Lambda))$ and $\Hom_{\Lambda_{\psi}}(\cpx{T},\mt_B(B))$ are exact, due to the isomorphism $_{\Lambda_{\psi}}\Lambda_{\psi}\simeq \mt_{\Lambda}(\Lambda)\oplus \mt_B(B)$.

To show that $\Hom_{\Lambda_{\psi}}(\cpx{T},\mt_{\Lambda}(\Lambda))$ is exact, we consider the exact sequence
$$0\lra \mz_{\Lambda}(I) \lraf{((0,1),0)} \mt_{\Lambda} (\Lambda) \lra \mt'_{\Lambda}(\Lambda) \lra 0$$
of $\Lambda_{\psi}$-modules, where $\mt'_{\Lambda}(\Lambda):=(\Lambda,M,\mu,0)$ and
$\mu: {}_BM\otimes_{A}\Lambda \ra {}_BM$ is the multiplication map.

Since $T^{i}$ is a projective $\Lambda_{\psi}$-module, the sequence of complexes
$$0\lra \Hom_{\Lambda_{\psi}}(\cpx{T},\mz_{\Lambda}(I)) \lra  \Hom_{\Lambda_{\psi}}(\cpx{T}, \mt_{\Lambda} (\Lambda)) \lra \Hom_{\Lambda_{\psi}}(\cpx{T},\mt'_{\Lambda}(\Lambda)) \lra 0$$
is exact. Thus, to show the exactness of $\Hom_{\Lambda_{\psi}}(\cpx{T}, \mt_{\Lambda} (\Lambda))$, it is sufficient to prove the one of the complexes $\Hom_{\Lambda_{\psi}}(\cpx{T},\mz_{\Lambda}(I))$ and $\Hom_{\Lambda_{\psi}}(\cpx{T}, \mt'_{\Lambda}(\Lambda))$. However, it follows from  Lemmas \ref{homo-equ}(2) and  \ref{homo-equ}(7) that $\Hom_{\Lambda_{\psi}}(\cpx{T},\mz_{\Lambda}(I))\simeq \Hom_{\Lambda}(\cpx{P},I)$, while the latter complex $\Hom_{\Lambda}(\cpx{P},I)$ is indeed exact, according to the weak compatibility of $_{\Lambda}I_{\Lambda}$.
Thus $\Hom_{\Lambda_{\psi}}(\cpx{T},\mz_{\Lambda}(I))$ is exact.
Further, it follows from the exact sequence
$$0\lra \mz_{B}(M) \lra \mt'_{\Lambda}(\Lambda) \lra \mz_{\Lambda}(\Lambda) \lra 0 $$
of $\Lambda_{\psi}$-modules and the projectivity of $T^i$ that the following sequence of complexes is exact
$$ 0\lra \Hom_{\Lambda_{\psi}}(\cpx{T}, \mz_{B}(M)) \lra
\Hom_{\Lambda_{\psi}}(\cpx{T},\mt'_{\Lambda}(\Lambda)) \lra
\Hom_{\Lambda_{\psi}}(\cpx{T}, \mz_{\Lambda}(\Lambda)) \lra 0.$$
By Lemmas \ref{homo-equ}(6)-(7),  $\Hom_{\Lambda_{\psi}}(\cpx{T}, \mz_{B}(M)) \simeq \Hom_{B}(\cpx{Q},M)$. Since ${}_{B}M_A$ is weakly compatible, the complex $\Hom_{B}(\cpx{Q},M)$ is exact, and therefore $\Hom_{\Lambda_{\psi}}(\cpx{T}, \mz_{B}(M))$ is exact. Similarly, by Lemmas \ref{homo-equ}(2) and \ref{homo-equ}(7), $ \Hom_{\Lambda_{\psi}}(\cpx{T}, \mz_{\Lambda}(\Lambda)) \simeq \Hom_{\Lambda}(\cpx{P}, \Lambda)$ is exact. Thus $\Hom_{\Lambda_{\psi}}(\cpx{T}, \mt'_{\Lambda}(\Lambda))$ is exact, and so $\Hom_{\Lambda_{\psi}}(\cpx{T},\mt_{\Lambda}(\Lambda))$ is exact.

Now, we show that $\Hom_{\Lambda_{\psi}}(\cpx{T}, \mt_{B}(B))$ is exact.
Similarly, from the exact sequence
$ 0\ra \mz_{\Lambda}(N) \ra \mt_{B}(B) \ra \mz_{B}(B) \ra 0 $
of $\Lambda_{\psi}$-modules we get the exact sequence of complexes
$$ 0\lra \Hom_{\Lambda_{\psi}}(\cpx{T}, \mz_{\Lambda}(N)) \lra
\Hom_{\Lambda_{\psi}}(\cpx{T},\mt_{B}(B)) \lra
\Hom_{\Lambda_{\psi}}(\cpx{T}, \mz_{B}(B)) \lra 0.$$
By Lemmas \ref{homo-equ}(2) and \ref{homo-equ}(7), together with the weak compatibility of $_{\Lambda}N_B$, we can show that the complex $\Hom_{\Lambda_{\psi}}(\cpx{T}, \mz_{\Lambda}(N))$ is exact. By Lemmas \ref{homo-equ}(6) and \ref{homo-equ}(7),
$\Hom_{\Lambda_{\psi}}(\cpx{T}, \mz_{B}(B)) \simeq \Hom_{B}(\cpx{Q},B)$. Since $\cpx{Q}$ is a totally exact complex,  $\Hom_{B}(\cpx{Q},B)$ is exact. Hence $\Hom_{\Lambda_{\psi}}(\cpx{T}, \mz_{B}(B)) $ is exact, and so is the complex $\Hom_{\Lambda_{\psi}}(\cpx{T},\mt_{B}(B))$. Thus the complex $\cpx{T}\in\mathscr{C}(\pmodcat{\Lambda_{\psi}})$ is totally exact.

Finally, we show  $\Ker(d^0_T)\simeq (X,Y,f,g)$. Clearly, $\Ker(d_F^0)=X$ and $\Ker(d^0_Y)=Y$.  Moreover, we can verify the following two exact commutative diagrams:
$$ \xymatrix{ & M\otimes_AX \ar[r]^{1_M\otimes d_X}\ar[d]^-{f} & M\otimes_AF^0\ar[r]^-{1_M\otimes d_F^0}\ar[d]^-{\big(\begin{smallmatrix}\pi_{P^0} & 0\\ 0 &0 \end{smallmatrix}\big)} & M\otimes_AF^1\ar[d]^-{\big(\begin{smallmatrix}\pi_{P^1} & 0\\ 0 &0 \end{smallmatrix}\big),} \\
0\ar[r] &Y\ar[r] & Y^0\ar[r]^-{d_Y^0} & Y^1
} \; \xymatrix{ & N\otimes_BY \ar[r]^{1_N\otimes d_Y}\ar[d]^-{g} & N\otimes_BY^0\ar[r]^-{1_N\otimes d_Y^0}\ar[d]^-{\big(\begin{smallmatrix}\psi_{P^0} & 0\\ 0 &1_{N\otimes_BQ^0} \end{smallmatrix}\big)} & N\otimes_BY^1\ar[d]^-{\big(\begin{smallmatrix}\psi_{P^1} & 0\\ 0 &1_{N\otimes_BQ^1} \end{smallmatrix}\big)} \\
0\ar[r] &X\ar[r]_-{d_X} & F^0\ar[r]_-{d_F^0} & F^1
}$$
This shows $\Ker(d^0_T) \simeq (X,Y,f,g)$. Thus $(X,Y,f,g)$ is a Gorenstein-projective $\Lambda_{\psi}$-module with a total projective resolution $\cpx{T}$. $\square$

\subsection{Necessary conditions for Gorenstein-projective modules\label{sect3.2}} 
In this section we discuss the converse of Theorem \ref{main-result-1}. We start with the following lemma

\begin{Lem} {\rm \cite[Proposition 6.1]{KT}}\label{tensor-equ}
 Let $U :=(C_A, D_B, h, k)$ be a right $\Lambda_{(\phi,\psi)}$-module, with $h\in \Hom_{B^{\opp}}(C\otimes_AN,D_B)$ and
 $k\in\Hom_{A^{\opp}}(D\otimes_BM_A, C_A)$,  and $V := (X, Y,f,g)$ a left $\Lambda_{(\phi,\psi)}$-module. Then there is an isomorphism of abelian groups $$ U \otimes_{\Lambda_{(\phi,\psi)}}V = (C\otimes_{A}X \oplus D\otimes_{B}Y) /H,$$
where $H$ is a subgroup of $C\otimes_{A}X \oplus D\otimes_{B}Y$ generated by
$\{c\otimes (n\otimes y)g - (c\otimes n)h \otimes y\mid c\in C, n\in N, y\in Y\}\cup \{ d \otimes (m\otimes x)f - (d\otimes m)k \otimes x\mid d\in D, x\in X, m\in M\}.$
\end{Lem}

\begin{Lem} \label{tensor}
Let $C\in \Lambda^{op}\modcat$, $X \in \Lambda \modcat$, $D\in B^{op}\modcat$, and $Y\in B\modcat$. Then there are the following isomorphisms of abelian groups, which are natural in each variables.

$(1)$ $\mz_{\Lambda^{op}}(C)\otimes_{\Lambda_{\psi}} \mt_{\Lambda}(X)\simeq C\otimes_{\Lambda}X$, $(c,0)\otimes ((x,i\otimes x'),m\otimes x'')\mapsto c\otimes x$, $x,x',x''\in X$, $c\in C$, $i\in I$, $m\in M$.
			
$(2)$ $\mz_{B^{op}}(D)\otimes_{\Lambda_{\psi}}   \mt_B(Y) \simeq D\otimes_BY$,  $(0,d)\otimes (n\otimes y',y)\mapsto d\otimes y$, $y,y'\in Y$, $d\in D$,  $n\in N$.
			
$(3)$ $\mz_{\Lambda^{op}}(C) \otimes_{\Lambda_{\psi}}  \mt_B(Y)=0$.
			
$(4)$ $\mz_{B^{op}}(D)\otimes_{\Lambda_{\psi}} \mt_\Lambda(X)=0$.
\end{Lem}

{\bf Proof.} We prove only (1) and (3), while the rest can be proved similarly and are left to the reader.

(1) Since $_AA\otimes_{\Lambda}X\simeq X(I)$ as $A$-modules, we get $C\otimes_{A}\, X(I)\simeq C\otimes_AA\otimes_{\Lambda}X\simeq C\otimes_{\Lambda}X$. By definition, $\mz_{\Lambda^{op}}(C)= (C_A,0,0,0)$ and $\mt_{\Lambda}(X)=(X(I), M\otimes_{\Lambda}X, \pi_X,\psi_X)$. By Lemma \ref{tensor-equ}, $\mz_{\Lambda^{op}}(C)\otimes_{\Lambda_{\psi}}\mt_{\Lambda}(X) \simeq (C\otimes_{\Lambda}X)/H$, while the subgroup $H$ is generated by $\{c\otimes (n\otimes m)\psi\; x)\mid c\in C, n\in N, m\in M, x\in X\}$. Thanks to $IX=0$, we get $H=0$. Thus (1) holds. Precisely, we can define $\alpha: C\otimes_{\Lambda}X \ra \mz_{\Lambda^{op}}(C)\otimes_{\Lambda_{\psi}} \mt_{\Lambda}(X)$ by $c\otimes x \mapsto (c,0)\otimes((x,0),0)$, and $\beta: \mz_{\Lambda^{op}}(C)\otimes_{\Lambda_{\psi}} \mt_{\Lambda}(X)\ra C\otimes_{\Lambda}X$, $(c,0)\otimes ((x,i\otimes x'),m\otimes x'')\mapsto c\otimes x$, $x,x',x''\in X$, for $c\in C$, $i\in I$, $m\in M$. One can check that they are homomorphisms of abelian groups, satisfying $\alpha\beta=1$ and $\beta\alpha=1$. Clearly, the isomorphisms of $\alpha$ and $\beta$ are natural in $C$ and $X$.

$(3)$ In this case, $H=C\otimes_{A}N\otimes_{B}Y$ in Lemma \ref{tensor-equ}. Thus $\mz_{\Lambda^{op}}(C) \otimes_{\Lambda_{\psi}}  \mt_B(Y) \simeq (C\otimes_{A}N\otimes_{B}Y)/H=0$. $\square$

\medskip
In the \textbf{rest of this section} we assume that $\cpx{T}=(T^i,d_T^i)_{i\in\mathbb{Z}} \in \mathscr{C}(\pmodcat{\Lambda_{\psi}})$ is a totally exact complex such that $\Ker(d_T^0)=(X,Y, f, g)\in \Lambda_{\psi} \modcat$.

By Lemma \ref{prop:projmod}, $T^{i}=\mt_{\Lambda}(P^i)\oplus \mt_{B}(Q^i)$ for some $P^i\in \pmodcat{\Lambda}$ and $Q^i\in
\pmodcat{B}$. Thus we may write precisely $d_T^i=(\begin{smallmatrix}
 t_{11}^{i} & t_{12}^{i}  \\
 t_{21}^{i} & t_{22}^{i}
\end{smallmatrix})$, with $t_{11}^{i} \in \Hom_{\Lambda_{\psi}}(\mt_{\Lambda}(P^{i}), \mt_{\Lambda}(P^{i+1}))$, $t_{12}^{i} \in \Hom_{\Lambda_{\psi}}(\mt_{\Lambda}(P^{i}), \mt_B(Q^{i+1}))$,
$t_{21}^{i}\in$ $\Hom_{\Lambda_{\psi}} (\mt_{B}(Q^{i}),$ $\mt_{\Lambda}(P^{i+1}))$, and $t_{22}^{i}\in \Hom_{\Lambda_{\psi}} (\mt_{B}(Q^{i}), \mt_B(Q^{i+1}))$.
By Lemma \ref{homo-equ}(3), there is a $\Lambda$-module homomorphism $(d_{P}^{i}, \alpha^{i})$, where $d_{P}^{i}: P^{i}\ra P^{i+1}$ and $\alpha^{i}: P^{i}\ra I\otimes_{\Lambda}P^{i+1}$ are homomorphisms of $\Lambda$-modules, such that $$t_{11}^{i}=\begin{pmatrix}
 d_{P}^{i} & \alpha^{i}  \\
 0 & 1_{I}\otimes d^i_P
\end{pmatrix},$$
Similarly, by Lemmas \ref{homo-equ}(1), \ref{homo-equ}(4) and \ref{homo-equ}(5),  we have a $\Lambda$-module homomorphism $\beta^{i}: P^{i}\ra N\otimes_{B}Q^{i+1} $, and two  $B$-homomorphisms $\rho^{i}: Q^{i}\ra M\otimes_{\Lambda}P^{i+1}$ and $d_{Q}^{i}: Q^{i}\ra Q^{i+1}$ of $B$-modules, such that
$t_{12}^{i}=\big(\binom{\beta^{i}}{0},0\big), \quad t_{21}^{i}=((0,\tau^{i}),\rho^{i}) \mbox{ with } \tau^i=(1_N\otimes \rho^{i})\psi_{{P^{i+1}}},$ and $t_{22}^{i}=(1_N\otimes d_{Q}^{i},d_{Q}^{i})$. Further, the exact complex $\cpx{T}$ provides two exact complexes $\cpx{F}:=(F^i,d_{F}^{i})\in \mathscr{C}(A\modcat)$ and $\cpx{Y}:=(Y^i,d_{Y}^{i})\in \mathscr{C}(B\modcat)$ by defining
$$(\ddag)\qquad
F^{i}:=P^i(I)\oplus (N\otimes_{B}Q^i), \; d_F^{i}:=\begin{pmatrix}
 (\begin{smallmatrix} d_{P}^{i} & \alpha^{i}       \\
 0 & 1_{I}\otimes d^i_P \end{smallmatrix}) & (\begin{smallmatrix}  \beta^{i} \\
 0\end{smallmatrix}) \\
  ( 0, \tau^i)& 1_{N}\otimes d^i_Q
\end{pmatrix},  \quad
Y^i:=M\otimes_\Lambda P^{i}\oplus Q^{i}, \; d_{Y}^{i}:=
\begin{pmatrix}
	 	1_M\otimes d_P^{i} & 0 \\
	 	\rho^{i} & d_ Q^{i}
	 \end{pmatrix}$$
Then  $\Ker(d_F^0)=X$ and $\Ker(d_Y^0)=Y$.

Let $\cpx{Q}:=(Q^i, d_Q^i)$ and $\cpx{P}: =(P^i, d_P^i)$.  Then it follows from $d_Y^id_Y^{i+1}=0$ and $d_F^id_F^{i+1}=0$ that $d_Q^id_Q^{i+1}=0$ and $d_P^id_P^{i+1}=0$, respectively. Thus $\cpx{Q}\in \C{\pmodcat {B}}$ and $\cpx{P}\in \C{\pmodcat {\Lambda}}$.  Moreover, we define $\cpx{Z}=(Z^i, d_Z^i)$ with
~$Z^{i}=(I\otimes_\Lambda P^{i})\oplus (N\otimes_BQ^{i})\in \Lambda\modcat$,
$d_{Z}^{i}=\bigl(\begin{smallmatrix}
 1_I \otimes d^i_P & 0 \\
 \tau^{i} & 1_N \otimes d_Q^{i}\end{smallmatrix}\bigr)$.
Then it follows again from $d_F^id_F^{i+1}=0$ that $d_Z^id_Z^{i+1}=0$. So  $\cpx{Z}\in \C \Lambda$.

Note that the complex $\cpx{Y}$ gives rise to an exact sequence in $\C{B\modcat}$:
\begin{equation} 0\lra M\otimes_{\Lambda} \cpx{P} \lraf{\cpx{a}} \cpx{Y} \lraf{\cpx{b}} \cpx{Q}\lra 0,  \label{ex1}  \end{equation}
where $\cpx{a}$ and $\cpx{b}$ are canonical inclusion and projection, respectively. Also, we have two exact sequences of complexes in $\C \Lambda$:
\begin{equation} 0\lra I\otimes_{\Lambda} \cpx{P} \lraf{\cpx{c}} \cpx{Z} \lraf{\cpx{d}} N\otimes_B \cpx{Q}\lra 0, \label{ex2}\end{equation}
\begin{equation} 0\lra  \cpx{Z} \lraf{\cpx{p}} \cpx{F} \lraf{\cpx{q}}  \cpx{P}\lra 0 \qquad \label{ex3}\end{equation}
where $\cpx{c}$ and $\cpx{p}$ are canonical inclusions, and where $\cpx{d}$ and $\cpx{q}$ are canonical projections.
Further, there is a chain map $\cpx{\sigma}$ in $\C{ \Lambda\modcat}$
$$\cpx{\sigma}=(\sigma^i)_{i\in \mathbb{Z}}: N\otimes_{B}\cpx{Y} \ra \cpx{Z},  \quad
\sigma^i:=\bigl(\begin{smallmatrix}
  \psi \otimes 1_{P^{i}} & 0 \\
 0 & 1_{N}\otimes 1_{Q^i}
\end{smallmatrix}\bigr),$$
such that the following diagram of complexes of $\Lambda$-modules is commutative and exact
\begin{equation} \label{ex4}
\xymatrix{
&N\otimes_{B}M\otimes_{\Lambda}\cpx{P} \ar[r]^-{1_N\otimes \cpx{a}} \ar[d]^{\psi\otimes 1_{\cpx{P}}} & N\otimes_{B}\cpx{Y} \ar[d]^{\cpx{\sigma}}\ar[r]^{1_N\otimes \cpx{b}} & N\otimes_{B}\cpx{Q} \ar[r] \ar@{=}[d]& 0\\
0\ar[r] & I\otimes_{\Lambda}\cpx{P}\ar[r]^{\cpx{c}} & \cpx{Z} \ar[r]^{\cpx{d}} &  N\otimes_{B}\cpx{Q} \ar[r] & 0 &&  }
\end{equation}

\begin{Lem} \label{toa-acy-deduce}
$(1)$  If $\mz_{\Lambda^{op}}(M)$ and $\mz_{\Lambda}(N)$ are semi-weakly compatible $\Lambda_{\psi}$-modules, then $\cpx{Q}\in \C{\pmodcat{B}}$  is totally exact.

$(2)$ If $\mz_{\Lambda^{op}}(I)$, $\mz_{\Lambda}(I)$, $\mz_{B^{op}}(N)$, and $\mz_{B}(M)$ semi-weakly compatible $\Lambda_{\psi}$-modules, then $\cpx{P}\in \C{\pmodcat{\Lambda}}$ is totally exact.
\end{Lem}
{\bf Proof.} $(1)$ Since $\cpx{T}\in \mathscr{C}(\pmodcat{\Lambda_{\psi}})$ is a totally exact complex and since $\mz_{\Lambda^{op}}{(M)}$ is semi-weakly compatible by assumption, the complex $\mz_{\Lambda^{op}}(M)\otimes_{\Lambda_{\psi}} \cpx{T}$ is exact. By Lemmas \ref{tensor}(1) and \ref{tensor}(3), $\mz_{\Lambda^{op}}(M)\otimes_{\Lambda_{\psi}} \cpx{T} \simeq M\otimes_{\Lambda}\cpx{P}$ as complexes, this yields that $M\otimes_{\Lambda}\cpx{P}$ is exact. It follows from the exact sequence (\ref{ex1}) that $\cpx{Q}$ is an exact complex. Since $\mz_{\Lambda}{(N)}$ is semi-weakly compatible, $ \Hom_{\Lambda_{\psi}}(\cpx{T}, \mz_{\Lambda}(N))$ is exact.  As $\cpx{T}$ is totally exact, the complex $\Hom_{\Lambda_{\psi}}(\cpx{T},\mt_{B}(B))$ is exact. Applying $\Hom_{\Lambda_{\psi}}(\cpx{T},-)$ to the following exact sequence
$ 0\ra \mz_{\Lambda}(N) \ra \mt_{B}(B) \ra \mz_{B}(B) \ra 0,$ we get the exact sequence of complexes of $\mathbb{Z}$-modules:
$$ 0\lra \Hom_{\Lambda_{\psi}}(\cpx{T}, \mz_{\Lambda}(N)) \lra
\Hom_{\Lambda_{\psi}}(\cpx{T},\mt_{B}(B)) \lra
\Hom_{\Lambda_{\psi}}(\cpx{T}, \mz_{B}(B)) \lra 0.$$
It follows that $\Hom_{\Lambda_{\psi}}(\cpx{T}, \mz_{B}(B))$ is exact. Similarly, by Lemmas \ref{homo-equ}(6)-(7),  the complex $\Hom_{B}(\cpx{Q},B) \simeq \Hom_{\Lambda_{\psi}}(\cpx{T}, \mz_{B}(B))$ is exact. Hence $\cpx{Q}$ is a totally exact complex.

(2) Since $\cpx{T}\in \mathscr{C}(\pmodcat{\Lambda_{\psi}})$ is totally exact, it follows from Lemmas \ref{tensor}(2) and \ref{tensor}(4), together with the semi-weak compatibility condition on $\mz_{B^{op}}(N)$, that $\mz_{B^{op}}(N)\otimes_{\Lambda_{\psi}} \cpx{T} \simeq N\otimes_{B}\cpx{Q}$ is exact. Similarly, by Lemma \ref{tensor}(1) and the assumption on $\mz_{\Lambda^{op}}(I)$, we know that $\mz_{\Lambda^{op}}(I)\otimes_{\Lambda_{\psi}} \cpx{T} \simeq I\otimes_{\Lambda}\cpx{P}$ is exact. It then follows from the exact sequence (\ref{ex2}) that $\cpx{Z}$ is exact. This implies, together with the exact sequence $(\ref{ex3})$, that $\cpx{P}$ is an exact complex.
By the semi-weak compatibility of $\mz_{\Lambda}(I)$,  we deduce that $\Hom_{\Lambda_{\psi}}(\cpx{T},\mz_{\Lambda}(I))$ is exact.
From the exact sequence $0\ra \mz_{\Lambda}(I) \lraf{((0,1),0)} \mt_{\Lambda} (\Lambda) \ra \mt'_{\Lambda}(\Lambda) \ra 0$ of $\Lambda_{\psi}$-modules, we have the exact sequence
$$ 0\lra \Hom_{\Lambda_{\psi}}(\cpx{T},\mz_{\Lambda}(I)) \lra  \Hom_{\Lambda_{\psi}}(\cpx{T}, \mt_{\Lambda} (\Lambda)) \lra \Hom_{\Lambda_{\psi}}(\cpx{T},\mt'_{\Lambda}(\Lambda)) \lra 0. $$
This shows that $\Hom_{\Lambda_{\psi}}(\cpx{T}, \mt'_{\Lambda}(\Lambda))$ is exact.
Now, applying $\Hom_{\Lambda_{\psi}}(\cpx{T},-)$ to the exact sequence
$0\ra \mz_{B}(M) \ra \mt'_{\Lambda}(\Lambda) \ra \mz_{\Lambda}(\Lambda) \ra 0, $ we  obtain the exact sequence
$$ 0\lra \Hom_{\Lambda_{\psi}}(\cpx{T}, \mz_{B}(M)) \lra
\Hom_{\Lambda_{\psi}}(\cpx{T},\mt'_{\Lambda}(\Lambda)) \lra
\Hom_{\Lambda_{\psi}}(\cpx{T}, \mz_{\Lambda}(\Lambda)) \lra 0.$$
As $\mz_{B}(M)$ is semi-weakly compatible, $\Hom_{\Lambda_{\psi}}(\cpx{T}, \mz_{B}(M))$ is exact, and therefore so is $\Hom_{\Lambda_{\psi}}(\cpx{T}, \mz_{\Lambda}(\Lambda))$. By Lemmas \ref{homo-equ}(2) and \ref{homo-equ}(7), $ \Hom_{\Lambda}(\cpx{P},\Lambda)\simeq \Hom_{\Lambda_{\psi}}(\cpx{T}, \mz_{\Lambda}(\Lambda))$ is exact. Hence $\cpx{P}$ is totally exact. $\square$

\begin{Lem} \label{Lambda-B-wcpx} Assume that $_{\Lambda}N_B$ and $_BM_{\Lambda}$ are weakly compatible. If 	$\mz_{\Lambda^{op}}(M)$ and $\mz_{\Lambda}(N)$ are semi-weakly compatible $\Lambda_{\psi}$-modules, then
$\mz_{B^{op}}(N)$ and $\mz_{B}(M)$ are semi-weakly compatible $\Lambda_{\psi}$-modules.	
\end{Lem}	

{\bf Proof.} Since $\cpx{T}$ is a totally exact complex of projective $\Lambda_{\psi}$-modules, we deduce from Lemma \ref{toa-acy-deduce}(1) that $\cpx{Q}$ is a totally exact complex. By assumption, $N_B$ is semi-weakly compatible. It follows from $\mz_{B^{op}}(N)\otimes_{\Lambda_{\psi}}\cpx{T} \simeq N\otimes_{B}\cpx{Q}$ that $\mz_{B^{op}}(N)\otimes_{\Lambda_{\psi}}\cpx{T}$ is exact.  Since $_BM$ is semi-weakly compatible and  $\Hom_{\Lambda_{\psi}}(\cpx{T},\mz_{B}(M)) \simeq \Hom_{B}(\cpx{Q},M)$, we see that $\Hom_{\Lambda_{\psi}}(\cpx{T},\mz_{B}(M))$ is exact. $\square$

\begin{Theo}\label{main-result-2}	
Assume that $_{\Lambda}N_B$, $_BM_{\Lambda}$ and $_{\Lambda}I_{\Lambda}$ are weakly compatible bimodules and that $\mz_{\Lambda^{op}}(M)$ and $\mz_{\Lambda}(N)$ are semi-weakly $\Lambda_{\psi}$-modules. Further, assume that $\mz_{\Lambda^{op}}(I)_{\Lambda_{\psi}}$ and  $_{\Lambda_{\psi}}\mz_{\Lambda}(I)$ are semi-weakly compatible. If a $\Lambda_{\psi}$-module$(X,Y,f,g)$ is Gorenstein-projective, then

$(a)$ $\Coker(f) \in B\Gpmodcat$ and $\Coker(g) \in \Lambda\Gpmodcat$, and

$(b)$ $\Img(f)\simeq M\otimes_{A} \Coker(g)$, $\Img(g) /IX\simeq N\otimes_{B} \Coker(f)$, and $IX\simeq I\otimes_{A} \Coker(g)$.
\end{Theo}

{\bf Proof.} Suppose $(X,Y,f,g) \in \Lambda_{\psi} \Gpmodcat$. Then there is a totally exact complex $\cpx{T}\in \mathscr{C}(\pmodcat{\Lambda_{\psi}})$, such that $\Ker(d_T^0) = (X,Y,f,g)$. By Lemmas \ref{Lambda-B-wcpx} and \ref{toa-acy-deduce}, the foregoing complexes $\cpx{P}\in \C{\pmodcat{\Lambda }}$ and $\cpx{Q}\in \C{\pmodcat{B}}$ are totally exact.
Thus $U:=\Ker (d_P^0)\in \Lambda \Gpmodcat$ and $V:=\Ker (d_Q^0)\in B\Gpmodcat$.

First, we show that $\Coker(g)$ is a Gorenstein-projective $\Lambda$-module.

Since $\Ker(d^0_F)=X$ (see notation in ($\ddag$)), we write $d_X: X\ra P^{0}(I) \oplus (N\otimes_{B}Q^{0})$ for the inclusion of $A$-modules. The restriction of $d_X$ to $\Lambda$-modules will be denoted by $d_{{}_{\Lambda}X}$. Then $d_{_{\Lambda}X}=(e_0,e_1,e_2): {}_{\Lambda}X \ra P^{0}\oplus ({}_{\Lambda}I\otimes_{\Lambda}P^{0}) \oplus ({}_{\Lambda}N\otimes_{B}Q^{0})$. We show that $e_i$ has the property:
$(ix)e_0=0$, $(ix)e_1=i\otimes(x)e_0$ and $(ix)e_2=0$ for $i\in I, x\in X$.

Indeed, the homomorphism $d_X$ of $A$-modules shows that $(((ax)e_0,(ax)e_1),(ax)e_2)$ = $a[(x)((e_0,e_1),e_2)]= a[((x)e_0,(x)e_1),(x)e_2]$ for $x\in X$, $a=(\lambda,i)\in A$. Further, the $\Lambda$-homomorphisms $e_i$ show the equality
$$ (((ax)e_0,(ax)e_1),(ax)e_2)=((\lambda(x)e_0+(ix)e_0,\lambda(x)e_1+(ix)e_1 ), \lambda(x)e_2+(ix)e_2). $$
By the $A$-module structure of $P^{0}(I)$ and $N\otimes_{B} Q^{0}$, one obtains immediately
$$ a[((x)e_0,(x)e_1),(x)e_2]=(\lambda,i)(((x)e_0,(x)e_1),(x)e_2)=((\lambda(x)e_0,   \lambda(x)e_1 +i\otimes(x)e_0),\lambda(x)e_2) ,$$
that is, $(ix)e_0=0$, $(ix)e_1=i\otimes(x)e_0$ and $(ix)e_2=0$.

Let $d_U$ be the inclusion of $_{\Lambda}U$ into $_{\Lambda}P^0$. It follows from the sequence (\ref{ex3}) that there is an exact commutative diagram of $\Lambda $-modules:
$$
\xymatrix{
0 \ar[r] & \Ker(d_Z^0) \ar[r]^-{\epsilon_{X}} \ar[d]^{ d_{Z}} & {}_{\Lambda}X \ar[d]^{(e_0,e_1,e_2)}\ar[r]^{\lambda_X} & U\ar[r] \ar[d]^{d_U}& 0  \\
0\ar[r] & I\otimes_{\Lambda}P^{0} \oplus N\otimes_{B} Q^{0} \ar[r]^-{\left(\begin{smallmatrix}
  0& 1 &0 \\    0 &0&1
\end{smallmatrix}\right)} &  P^{0} \oplus I\otimes_{\Lambda}P^{0} \oplus N\otimes_{B} Q^{0} \ar[r]^-{\left(\begin{smallmatrix}
   1 \\   0 \\   0
\end{smallmatrix}\right)} &  P^{0} \ar[r] & 0 }
$$
where all vertical maps are injective and $e_0=\lambda_{X}d_U$. Note that $(ix)\lambda_{X}d_U=(ix)e_0=0$ for $i\in I$ and $x\in X$. Since $d_U$ is injective, one must have $(ix)\lambda_{X}=0$. Now, if we consider $_{\Lambda}U$ as an $A$-module, that is $IU=0$, then $\lambda_{X}$ is a homomorphism of $A$-modules.

Furthermore, $\epsilon_X$ is a homomorphism of $A$-modules if $_{\Lambda}\Ker(d^0_Z)$ is regarded as an $A$-module.  In fact, for $z\in \Ker (d_Z^0)$, let $x_z:=(z)\epsilon_{X}$. Then it follows from $i\otimes (x_z)e_0=i\otimes (x_z)\lambda_{X}d_U=0$ that $(ix_z)(e_0,e_1,e_2)=(0, i\otimes (x_z)e_0,0)=0$. Since the map $(e_0,e_1,e_2)$ is injective, we obtain $ix_z=0$, that is, $I\Img(\epsilon_X)=0$. This implies that $\epsilon_X$ is a homomorphism of $A$-modules. Thus there is an exact sequence of $A$-modules
\begin{equation}
\xymatrix{0 \ar[r] & \Ker(d_Z^0) \ar[r]^-{\epsilon_{X}}  & X \ar[r]^{\lambda_X} & U\ar[r] & 0}                              \qquad \end{equation}
which fits into the exact commutative diagram of $A$-modules:
$$\xymatrix@C=0.5cm@R=0.5cm{
   & 0 \ar[d] & 0 \ar[d] & 0 \ar[d] & \\
  0 \ar[r] &  \Ker(d_Z^0) \ar[d]_{\epsilon_X} \ar[r]^-{d_Z} & (I\otimes_{\Lambda}P^{0})\oplus (N\otimes_{B}Q^{0}) \ar[d]^{j^0}
  \ar[r]^{d_Z^0} & (I\otimes_{\Lambda}P^{1})\oplus (N\otimes_{B}Q^{1}) \ar[d]^{j^1}    \\
 0 \ar[r] &  X \ar[r]^-{d_X} \ar[d]_{\lambda_{X}} & P^0(I)\oplus (N\otimes_{B}Q^{0}) \ar[r]^{d_E^0} \ar[d]^{k^0} &  P^1(I) \oplus (N\otimes_{B}Q^{1})  \ar[d]^{k^1}   \\
 0 \ar[r] &  U \ar[d] \ar[r]^-{d_U} & P^0 \ar[d] \ar[r]^{d^0_P} & P^1 \ar[d]   \\
  & 0 & 0 & 0 &  }
$$
with $j^i=\left(\begin{smallmatrix}
   (0,1)& 0 \\  0 & 1_{N\otimes_BQ^i}
\end{smallmatrix}\right)$ in which $(0,1):I\otimes_{\Lambda}P^{i} \ra P^i \oplus  I\otimes_{\Lambda}P^{i} $ is the canonical inclusion; and with $k^i=\left(\begin{smallmatrix}
    \binom{1}{0}\\   0
\end{smallmatrix}\right)$ in which $ \binom{1}{0}$ means the projection map
$P^i \oplus  I\otimes_{\Lambda}P^{i} \ra P^i$. Note that all homomorphisms of $\Lambda$-modules in the top and bottom rows are regarded as homomorphisms of $A$-modules via the canonical map $A\to \Lambda$. According to Lemmas \ref{A-Lambda-homo}(1)-(2), the vertical maps $j^i$ and $k^i$ are also homomorphisms of $A$-modules.

Since $\Ker(d_Y^0)=Y$, we have an inclusion  $d_Y: Y\hookrightarrow M\otimes_{\Lambda}P^0\oplus Q^0$. Let $\delta: N\otimes_{B}Y \ra \Ker(d_Z^0)$ be the homomorphism of $A$-modules induced from the chain map $\cpx{\sigma}$ in the diagram (\ref{ex4}).
Then $\delta$ is surjective. Actually, since $I_{\Lambda}$ is semi-weakly compatible and since we have shown that $\cpx{P}$ is a totally exact complex of projective modules, the complex $I\otimes_{\Lambda}\cpx{P}$ is exact. Similarly, $N\otimes_{B}\cpx{Q}$ is exact. So the exact sequence (\ref{ex3}) implies that $\cpx{Z}$ is exact. By the diagram (\ref{ex4}), the following diagram of $A$-modules is exact and commutative
$$
\xymatrix{
&N\otimes_{B}M\otimes_{\Lambda}U \ar[r]^-{1_N\otimes \eta_Y} \ar[d]^{\psi\otimes 1_U} & N\otimes_{B} Y \ar[d]^{\delta}\ar[r]^{1_N\otimes \mu_Y} & N\otimes_{B} V \ar[r] \ar@{=}[d]& 0\\
0\ar[r] & I\otimes_{\Lambda}U \ar[r] & \Ker(d_{Z}^0) \ar[r] &  N\otimes_{B} V \ar[r] & 0 }
$$
Therefore the Snake Lemma shows that $\delta$ is surjective.

Moreover, the diagram (\ref{ex4}) gives rise to the following one of $A$-modules
$$\xymatrix@C=0.8cm@R=1.1cm{
	& N\otimes_{B} Y \ar[dl]_{\delta} \ar[r]^-{1_N\otimes d_Y} \ar'[d]_(0.53){g}[dd]
	&   	(N\otimes M\otimes P^0)\oplus (N\otimes_BQ^0) \ar'[d]^{\bigl(\begin{smallmatrix}
			\psi_{{P^{0}}} & 0 \\
			0 & 1_{N\otimes{Q^0}}
		\end{smallmatrix}\bigr)}[dd]  \ar[r]^-{1_N\otimes d_Y^0}\ar[dl]_{\sigma^0} &  (N\otimes M\otimes P^1)\oplus (N\otimes_BQ^1)\ar[dl]^{\sigma^1} \ar[dd]^{\bigl(\begin{smallmatrix}
		\psi_{{P^{1}}} & 0 \\
		0 & 1_{N\otimes{Q^1}}
	\end{smallmatrix}\bigr)}     \\
	\Ker (d_Z^0) \ar[r]\ar[dr]_(0.5){\epsilon_{X}}
	&   (I\otimes P^0) \oplus (N\otimes_BQ^0)  \ar[dr]^(0.5){j^0} \ar[r]^(0.5){d_Z^0}  &(I\otimes P^1) \oplus (N\otimes_BQ^1) \ar[dr]^(0.5){j^1} \\
	& X \ar[r]_-{d_X}
	&   P^0(I) \oplus (N\otimes_BQ^0)    \ar[r]_{d_F^0}  &   P^1(I) \oplus (N\otimes_BQ^1)         }
$$
By the definitions of $g, \delta$ and $\epsilon_{X}$, the two top and two bottom squares are commutative. We can verify $\sigma^i j^i =\bigl(\begin{smallmatrix}
\psi_{{P^{i}}} & 0 \\ 0 & 1_{N\otimes{Q^i}}	\end{smallmatrix}\bigr)$ for all $i$.
Since $d_X$ is injective, there holds $\delta\,\epsilon_X=g$. Thus $\Coker (g)=\Coker(\epsilon_X) \simeq U\in \Lambda \Gpmodcat$.

Next, we prove $\Coker(f)\in B \Gpmodcat$. Observe that $M\otimes_{A}\cpx{P} \simeq M\otimes_{\Lambda}\cpx{P}$ as complexes of $B$-modules.  So the exact sequence (\ref{ex1}) may be rewritten as the following exact sequence of $B$-modules
$$0\lra M\otimes_{A} \cpx{P} \lraf{\cpx{a'}} \cpx{Y} \lraf{\cpx{b}} \cpx{Q}\lra 0$$
with $a'^i:  M\otimes_A P^i \lraf{\simeq}  M\otimes_{\Lambda} P^i \lraf{a^i} M\otimes_{\Lambda} P^i \oplus Q^i$.
This gives rise to the exact sequence of $B$-modules
$$0\lra M\otimes_{A}U \lraf{\eta'_{Y}} Y \lraf{\mu_Y} V\lra 0.$$
Now consider the diagram of $B$-modules
$$\xymatrix{
	& M\otimes_{A} X \ar[dl]_{1_{M}\otimes \lambda_{X}} \ar[r] \ar'[d]^(0.53){f}[dd]
	&   	( M\otimes_{A} P^0(I))\oplus (M\otimes_{A}N\otimes_BQ^0) \ar'[d]^{\bigl(\begin{smallmatrix}
			\pi_{{P^{0}}} & 0 \\
			0 & 0
		\end{smallmatrix}\bigr)}[dd]  \ar[r]^-{1_M\otimes d_F^0}\ar[dl]^{1_M \otimes k^0} &  ( M\otimes_{A} P^1(I))\oplus (M\otimes_{A}N\otimes_BQ^1)\ar[dl]^{1_M \otimes k^1} \ar[dd]^{\bigl(\begin{smallmatrix}
		\pi_{{P^{1}}} & 0 \\
		0 & 0
	\end{smallmatrix}\bigr)}     \\
	 M\otimes_{A}U \ar[r]\ar[dr]_(0.5){\eta'_{Y}}
	&   M\otimes_{A} P^0 \ar[dr]_(0.5){a'^0} \ar[r]_(0.5){1_M \otimes d_P^0}  &M\otimes_A P^1  \ar[dr]^(0.5){a'^1} \\
	& Y \ar[r]_{d_Y}
	&  M\otimes_{\Lambda} P^0 \oplus Q^0 \ar[r]_{d_Y^0}  &   M\otimes_{\Lambda} P^1 \oplus Q^1         }
$$
By the definition of $\lambda_X$ as an $A$-module homomorphism, the upper two squares are commutative. By the definition of $\eta'_{Y}$, the lower two squares are commutative. Moreover, $(1_M \otimes k^i)a^i=\bigl(\begin{smallmatrix}
\pi_{{P^i}} & 0 \\ 0 & 0	\end{smallmatrix}\bigr)$. It follows from the injective map $d_Y$ that $(1_M \otimes \lambda_X)\eta'_{Y}=f$. Therefore $\Coker(f)=\Coker(\eta'_{Y})=V \in B\Gpmodcat$ and $\Img(f)\simeq M\otimes_{A}U = M\otimes_{A} \Coker(g)$. This completes the proof of $(a)$.

Having proved that $\Img(f)\simeq M\otimes_{A} \Coker(g)$, we now prove $\Img(g)/IX\simeq N\otimes_B\Coker(f)$.

Recall that $\cpx{F}\in \C{A}$ stands for the complex defined  in ($\ddag$). Let $W^i:= {}_{\Lambda}P^i\oplus {}_{\Lambda}N \otimes_{B}Q^i$, $d_W^i:=\bigl(\begin{smallmatrix}
		d_P^{i} & \beta^i \\
		0 & 1_N \otimes d_Q^i
	\end{smallmatrix}\bigr)$.
Due to $d_F^{i}d_F^{i+1}=0$, we have $d_P^i\beta^{i+1}+\beta^i(1_N \otimes d_Q^i)=0$ and $d_W^id_W^{i+1}=0$. So $\cpx{W} \in \mathscr{C}(\Lambda)$. Regarding $\Lambda$-modules as $A$-modules, we have the exact sequence
$$ (\#)\qquad 0\lra N\otimes_B \cpx{Q}\lraf{\cpx{s}} \cpx{W} \lraf{\cpx{t}} \cpx{P} \lra 0$$
of complexes in $\C A$. Since $_{\Lambda}N_B$ is weakly compatible and $\cpx{Q}$ is a totaly exact complex in $\C{\pmodcat{B}}$, the complex $N\otimes_{B}\cpx{Q}$ is exact. It then follows from the exactness of $\cpx{P}$ that the complex $\cpx{W}$ is exact.
Now, since $_{A}\Lambda\otimes_{A}F^{i} = {}_A\Lambda\otimes_{A}(P^{i}(I)\oplus N\otimes_{B}Q^i)\simeq P^i \oplus (N\otimes_{B}Q^i)$ and $1_\Lambda \otimes_{A} d_{F}^{i}=d_W^i $, we have $\cpx{W}\simeq \Lambda\otimes_{A}\cpx{F}$ as complexes in $\C A$. Hence $_A\Lambda\otimes_{A}\cpx{F}$ is an exact complex of $A$-modules and $_{A}\Ker(d_W^0) \simeq {}_{A}\Lambda\otimes_{A}X\simeq X/IX$. Thus ($\#$) induces the exact sequence of $A$-modules:
$ 0\ra {}_AN\otimes_B V\lraf{s} X/IX \lraf{t} {}_AU\ra 0. $
It follows from $\Lambda\otimes_{A}\cpx{F}\simeq \cpx{W}$ that there is a canonical chain map $\cpx{p_F}: \cpx{F} \ra \cpx{W}$ in $\C A$, with $p_F^i=\begin{pmatrix}
\binom{1}{0} & 0 \\  0 & 1_N \otimes 1_{Q^i}
\end{pmatrix}$. Now, we take the kernels of $\cpx{p_F}$ at degree $0$  and  get the canonical projection $p_X:X\ra X/IX$.
Considering the commutative diagram of complexes in $\C A$:
$$\xymatrix{
N\otimes_{B}\cpx{Y} \ar[rr]^{1_N\otimes\cpx{b}} \ar[d]_{\bigl(\begin{smallmatrix}
		\psi_{P^{i}} & 0 \\
		0 & 1_N \otimes 1_{Q^i}
	\end{smallmatrix}\bigr)} &&N\otimes_{B}\cpx{Q}\ar[d]_{(0,1)}\\
\cpx{F} \ar[rr]^{\cpx{p_F}} &&\cpx{W}}$$
and the differentials in degree $0$ in the diagram, we get the exact commutative diagram
$$\xymatrix@C=0.8cm{
	& N\otimes_B Y \ar[dl]_{g}\ar[rr] \ar'[d]^(0.53){1_N\otimes_B \mu_Y}[dd]
	&  & N\otimes_B Y^0 \ar'[d]^{b^0}[dd]  \ar[rr]^(0.3){1_N\otimes_B d_Y^0}\ar[dl]_(0.66){\bigl(\begin{smallmatrix}
			\psi_{P^{i}} & 0 \\
			0 & 1_N \otimes 1_{Q^0}
		\end{smallmatrix}\bigr)} & & N\otimes_B Y^1 \ar[dd]^{b^1}  \ar[dl]_(0.66){\bigl(\begin{smallmatrix}
		\psi_{P^{i}} & 0 \\
		0 & 1_N \otimes 1_{Q^0}
	\end{smallmatrix}\bigr)}   \\
	X \ar[rr]\ar[dd]^(0.23){p_X}
	&  &  P^0(I) \oplus (N\otimes_B Q^0) \ar[dd]^(0.23){p_F^0} \ar[rr]_(0.6){d_F^0} & &P^1(I) \oplus (N\otimes_B Q^1) \ar[dd]^(0.23){p_F^1} \\
	&  N\otimes_B V \ar[dl]^{s}\ar'[r][rr]
	&  &    N\otimes_B Q^0  \ar[dl]^{(0,1)}\ar'[r]^{1_N\otimes_B d_Q^0}[rr] & &    N\otimes_B Q^1   \ar[dl]^{(0,1)}      \\
	X/IX \ar[rr]^{d_W}
	&  & P^0 \oplus (N\otimes_B Q^0) \ar[rr]^{d_W^0} & &  P^1 \oplus (N\otimes_B Q^1)  }
$$in $\C A$ by Lemma \ref{ker-seq-sur}(2).
Thus $gp_X=(1_N\otimes \mu_Y)s$. Since $s$ is injective and $1_N \otimes\mu_Y$ is surjective, we obtain $\Img(g)/IX\simeq N \otimes_{B}V=N \otimes_{B}\Coker(f)$.

Finally, we show $IX\simeq I\otimes_{A}\Coker(g)$. Actually, it follows from the chain map $\cpx{c}: I\otimes_{\Lambda}\cpx{P}\to\cpx{Z}$ in (\ref{ex4}) that the following exact commutative diagram exists
$$\xymatrix@C=0.5cm@R=0.5cm{
 & I\otimes_{\Lambda}U \ar[r]\ar[d]^-{c} & I\otimes_{\Lambda} P^0 \ar[r]\ar[d]^-{c^0} & I\otimes_{\Lambda}P^1\ar[d]^-{c^1}\\
0\ar[r] & \Ker(d_Z^0) \ar[r] & Z^0\ar[r] & Z^1\\
} $$

Now, consider the diagram
\begin{equation}\label{ex12}
\xymatrix@R=0.6cm{
N\otimes_{B}M\otimes_{A}X \ar@{->>}[d]_{1_{N\otimes_{B}M}\otimes \lambda_{X}}\ar[rr]^-{\psi\otimes 1_X} & & I\otimes_{A}X \ar[dl]^{1_I\otimes\lambda_X} \ar[ddd]^{\mlt_X}\\
N\otimes_{B}M\otimes_{A}U \ar[d]_-{\simeq} \ar@{->>}[r]^-{\psi\otimes1_U}& I\otimes_{A}U \ar[d]^{\simeq} &  \\
N\otimes_{B}M\otimes_{\Lambda}U \ar[d]_-{1_N\otimes \eta_Y} \ar@{->>}[r]^-{\psi\otimes1_U}& I\otimes_{\Lambda}U \ar[d]^{c} &  \\
N\otimes_{B}Y \ar@{->>}[r]^{\delta} & \Ker (d_Z^0) \ar[r]^{\epsilon_X} &X  &&  }\end{equation}
Note that the out-side square is  commutative, due to $(X,Y,f,g)\in \Lambda_{\psi}\Modcat$. The down-left square commutes because of the commutative diagram (\ref{ex4}), while the upper-left square commutes, due to the property of $\psi$. Thus it follows from the surjective map $\psi\otimes 1_X$ that the right-side square is commutative. Since $c$ and $\varepsilon_x$ are injective maps, we have $IX\simeq I\otimes_{A} U \simeq I\otimes_{A} \Coker(g)$. This completes the proof of $(b)$. $\square$

\begin{Theo} \label{main-result-3} The following are equivalent for the Morita context ring $\Lambda_{\psi}$.

$(1)$ ${}_\Lambda N_B$, $_BM_\Lambda$ and $_{\Lambda}I_{\Lambda}$ are weakly compatible bimodules, the left $\Lambda_{\psi}$-modules $(_AN,0,0,0)$ and $(_AI,0,0,0)$ and the right $\Lambda_{\psi}$-modules $(M_A, 0,0,0)$ and $(I_A, 0,0,0)$ are semi-weakly compatible.

$(2)$ A $\Lambda_{\psi}$-module $(X,Y,f,g)$ is Gorensteion-projective if and only if

$\quad (a)$  $_B\Coker(f)$ and $_{\Lambda}\Coker(g)$ are Gorenstein-projective, and	

$\quad (b)$	$_B\Img(f)\simeq {}_BM\otimes_A\Coker(g)$, $_A\Img(g)/IX\simeq {}_AN\otimes_B\Coker(f)$, and
$_AIX\simeq {}_AI\otimes_A\Coker(g)$.
\end{Theo}

{\bf Proof.} (1) $\Rightarrow$ (2).  This follows from Theorems \ref{main-result-1} and \ref{main-result-2}.

(2) $\Rightarrow$ (1). This will be done in the rest of this section. So, in the following we always assume (2).

\begin{Lem} \label{sublem1} If $G\in \Gpmodcat{\Lambda}$, then $\mt_{\Lambda}(G)\in \Gpmodcat{\Lambda_{\psi}}.$ Similarly, if $Q\in \Gpmodcat{B}$, then $\mt_{B}(Q)\in \Gpmodcat{\Lambda_{\psi}}.$
\end{Lem}

{\bf Proof.} Since $\mt_{\Lambda}(G)= (A\otimes_{\Lambda}G, M\otimes_{\Lambda}G, \pi_G, \psi_P)$ where $\pi_G: {}_BM\otimes_AA\otimes_{\Lambda}G\simeq {}_BM\otimes_{\Lambda}G$ and $\psi_G: {}_AN\otimes_B(M\otimes_{\Lambda}G)\ra {}_AA\otimes_{\Lambda}G$ is given by $ n\otimes(m\otimes x)\mapsto (n\otimes m)\psi \otimes x
 $ for $n\in N, m\in M, x\in G$. It follows from $\Coker(\pi_G)=0$ and $_A\Img(\psi_G) = I\otimes_{\Lambda}G$ that $_{\Lambda}\Coker(\psi_G)\simeq {}_{\Lambda}G$ and the condition (2)(a) is satisfied. On the other hand, $\Img(\pi_G)= M\otimes_{\Lambda}G\simeq M\otimes_{\Lambda}\Coker(\psi_G)\simeq M\otimes_{A}\Coker(\psi_G)$. We can show $I\otimes_AG= I(A\otimes_{\Lambda}G)$ in $A\otimes_{\Lambda}G$, and therefore $\Img(\psi_G)/I(A\otimes_{\Lambda}G)=0$
and $I\otimes_{A}\Coker(\psi_G) = I\otimes_A G\simeq  I(A\otimes_{\Lambda}G)$. This means that the condition (2)(b) is satisfied. Hence $\mt_{\Lambda}(G)$ is a Gorenstein-projective $\Lambda_{\psi}$-module by (2).

Let $_BQ$ be a Gorenstein-projective module. By definition, $\mt_{B}(Q)=(N\otimes_BQ,Q,\Phi_Q,id_{N\otimes Q}),$ we show that the conditions (2)(a)-(b) hold for $\mt_{B}(Q)$. Note that $\Phi_Q: M\otimes_AN\otimes_BQ\ra {}_BQ=0$. It is now easy to verify that (a) and (b) holds true. Thus $\mt_{B}(Q)$ is Gorenstein-projective by (2).
$\square$

\begin{Lem} \label{sublem2} $(i)$ If $G\in \Gpmodcat{\Lambda}$, then $\Tor^{\Lambda}_{1}(I,G)=\Tor^{\Lambda}_{1}(M,G)=0.$ Thus $\Tor^{\Lambda}_i(I,G)=\Tor^{\Lambda}_i(M,G)=0$ for all $i>0$.

$(ii)$ If $_BW\in B\Gpmodcat$, then $\Tor_{1}^{B}(N,W)=0$. Thus $\Tor^{\Lambda}_{i}(N,W)=0$ for all $i>0$.
\end{Lem}

{\bf Proof.} Assume that $\cpx{P}=(P^i, d^i_P)$ is a total projective resolution of $G$ with $\Ker(d^0_P)=G$. Let $H: = {}_{\Lambda}\Ker(d^{-1}_P)$ and $b': H\hookrightarrow P^{-1}$ be the inclusion. By Lemma \ref{sublem1}, $\mt_{\Lambda}(H)= (H(I), M\otimes_{\Lambda}H, \pi_H,\psi_H)\in \Gpmodcat{\Lambda_{\psi}}$. Since $H(I)$ is a finitely generated $A$-module,  $\Hom_A(H(I),{}_AA_A)$ is a finitely generated right $A$-module. Suppose that $f_1,f_2, \cdots, f_s$ form a set of generators for $\Hom_A(H(I),{}_AA_A)$. Then $\alpha: H(I)\lraf{(f_1,\cdots, f_s)} ({}_AA)^s$ is a left $\add(_AA)$-approximation of $H(I)$. So we assume that $\alpha: H(I)\ra Q$ be a left $\add(_AA)$-approximation of $H(I)$ with $Q\in\add(_AA)$. Since $A$ is the trivial extension of $\Lambda$ by $I$, we may further assume $Q=P(I)$ for some $P\in\add(_{\Lambda}\Lambda)$. Then, by Lemma \ref{homo-equ}(3), $\alpha$ is of the form $\alpha=\left(\begin{smallmatrix}a & d  \\ 0 & 1_I\otimes a
\end{smallmatrix}\right): H\oplus I\otimes_{\Lambda}H\lra P\oplus I\otimes_{\Lambda}P$, where $a: H\to P$ and $d: H\to I\otimes_{\Lambda}P$ are homomorphisms of $\Lambda$-modules. We show that $a:H\to P$ is an injective left $\add(_\Lambda\Lambda)$-approximation of $H$. Actually, for $b: H\to P'$ with $P'$ in $\pmodcat{\Lambda}$, we have an $A$-module homomorphism $\bar{b}:=\left(\begin{smallmatrix}b & 0  \\ 0 & 1_I\otimes b\end{smallmatrix}\right): H(I)\ra P'(I)$. Since $\alpha$ is an approximation, there is a homomorphism $\bar{c}=\left(\begin{smallmatrix}c & e  \\ 0 & 1_I\otimes c
\end{smallmatrix}\right):P(I)\ra P'(I)$ such that $\bar{b}=\alpha\bar{c}$, that is, $\left(\begin{smallmatrix}b & 0  \\ 0 & 1_I\otimes b
\end{smallmatrix}\right)=\left(\begin{smallmatrix}a & d  \\ 0 & 1_I\otimes a
\end{smallmatrix}\right)\left(\begin{smallmatrix}c & e  \\ 0 & 1_I\otimes c
\end{smallmatrix}\right)$. This implies that $b=ac$ and $a$ is a left $\add(_{\Lambda}\Lambda)$-approximation of $H$. Taking $b=b'$ , we see immediately that $a$ is injective.

Since $a$ is a left $\add(_{\Lambda}\Lambda)$-approximation of $H$, there is a homomorphism $c: P\to P^{-1}$ of $\Lambda$-modules such that the following diagram is exact and commutative
$$\xymatrix@C=0.4cm@R=0.7cm{
 0\ar[r]& H\ar[r]^-{a}\ar@{=}[d] & P\ar[r]\ar[d]^-{c} &\Coker(a)\ar[r]\ar@{-->}[d] & 0\\
0\ar[r] & H \ar@{^(->}[r] & P^{-1}\ar[r] & G\ar[r] &0.
} $$
Thus the right-side square is a pushout and pullback diagram. This induces an exact sequence of $\Lambda$-modules:
$$\xymatrix{ 0\ar[r] & P\ar[r]\ar[r] & P^{-1}\oplus \Coker(a)\ar[r] & G\ar[r] & 0
}$$
As ${}_{\Lambda}G$ is Gorenstein-projective, we always have $\Ext^i_{\Lambda}(G, X)=0$ for $i>0$ and any module $_{\Lambda}X$ of finite projective dimension. Thus the exact sequence splits and $\Coker(a)\oplus P^{-1}\simeq G\oplus P$ as $\Lambda$-modules.
Therefore, to show Lemma \ref{sublem2}, it sufficient to show that $\Tor_1^{\Lambda}(I, \Coker(a))=\Tor_1^{\Lambda}(M,\Coker(a))=0,$ that is,
the maps $1_I\otimes a$ and $1_M\otimes a$ are injective. This is done by considering the left $\add(_{\Lambda_{\psi}}\Lambda_{\psi})$-approximation of $\mt_{\Lambda}(H)$.

In fact, we have $_{\Lambda_{\psi}}\Lambda_{\psi}=\mt_A(A)\oplus \mt_B(B)$. Since $\mt_{A}$ is a fully faithful additive functor (see Lemma \ref{mr-recol}), it follows from the left $\add(_AA)$-approximation $\alpha: H(I) \ra P(I)$ of $H(I)$ that $\mt_{A}(\alpha): \mt_{A}(H(I)) \ra \mt_{A}(P(I))$ is a left $\add(\mt_{A}(A))$-approximation of $\mt_{A}(H(I))$. This implies also that $(\left(\begin{smallmatrix}a & d  \\ 0 & 1_I\otimes a
\end{smallmatrix}\right),1_M\otimes a):\mt_{\Lambda}(H)\ra \mt_{\Lambda}(P)$ is  a  left $\add(\mt_{\Lambda}(\Lambda))$-approximation of $\mt_{\Lambda}(H)$. Let $\beta= \left(\begin{smallmatrix}a & d  \\ 0 & 1_I\otimes a
\end{smallmatrix}\right),1_M\otimes a)$. Take a left $\add(\mt_{B}(B))$-approximation of $\mt_{\Lambda}(H)$, say $\theta: \mt_{\Lambda}(H) \ra \mt_{B}(Q)$ for some $Q\in \add(_BB)$. By Lemma \ref{homo-equ}(1), $\theta$ is of the form $(\binom{h}{0},0)$ with $h\in \Hom_{\Lambda}(H,N\otimes_{B}Q)$. Then we get a left $\add(_{\Lambda_{\psi}}\Lambda_{\psi})$-approximation $(\beta, \theta): \mt_{\Lambda}(H) \ra \mt_{\Lambda}(P) \oplus \mt_{B}(Q)$ of $\mt_{\Lambda}(H)$.
Since $\mt_{\Lambda}(H)\in \Lambda_{\psi}\Gpmodcat$ by Lemma \ref{sublem1}, there is an injective homomorphism from $\mt_{\Lambda}(H)$ to a projective $\Lambda_{\psi}$-module, and therefore $(\beta, \theta)$ is injective. This shows that the homomorphisms $\left(\begin{smallmatrix}a & d & h \\ 0 & 1_I\otimes a & 0
\end{smallmatrix}\right): H \oplus I\otimes_{\Lambda}H  \ra P \oplus I\otimes_{\Lambda}P  \oplus N\otimes_{B}Q $ of $\Lambda$-modules and $(1_M\otimes a,0): M\otimes_{\Lambda}H \ra M\otimes_{\Lambda}P\oplus Q$ of $B$-modules are injective. Therefore $1_I\otimes a$ and $1_M\otimes a$ are injective. Thus $\Tor_{1}^{\Lambda}(M,\Coker (a))=0$ and $\Tor_{1}^{\Lambda}(I,\Coker (a))=0$. Therefore $\Tor_{1}^{\Lambda}(I,G)=0$. A dimension shift argument shows $\Tor_{i}^{\Lambda}(I,G)=0$ for $i>0$.

Now, let $W\in B\Gpmodcat$ with $\cpx{Q}$ a totally exact complex in $\mathscr{C}(\pmodcat{B})$, such that $\Ker( d_Q^0) = W$. Then there is the short exact sequence
$  0\ra V \lraf{v} Q^{-1} \lraf{w} W \ra 0$
 of $B$-modules with $V=\Ker(d_Q^{-1})$. This yields an exact sequence of $\Lambda_{\psi}$-modules
$$    0\lra (U,V, s,t) \lraf{(i,v)}  \mt_{B}(Q^{-1}) \lraf{(1_N\otimes w,w)} \mt_{B}(W) \lra 0 $$
with $i: {}_AU\ra {}_AN\otimes_{B}Q^{-1}$ the kernel of $1_N\otimes w$.
By the diagram (\ref{morphismdiagrams}) (see Section \ref{sect2.1}), the homomorphisms $s$ and $t$ fit in the exact commutative diagrams, respectively:
$$\xymatrix{
   & M\otimes_A U \ar[r]^-{1_{M}\otimes i} \ar[d]^{ s} & M\otimes_AN\otimes_BQ^{-1} \ar[d]^{0}     \\
  0\ar[r] & V    \ar[r]^{ v} & \; Q^{-1},}\qquad  \xymatrix{
  & N\otimes_B V \ar[r]^{1_{N}\otimes v} \ar[d]^{  t}
  & N\otimes_B  Q^{-1} \ar@{=}[d] \ar[r]^{1_N \otimes w}
  & N\otimes_B W \ar@{=}[d] \ar[r] & 0\\
0\ar[r]&U\ar[r]^-{i} &  N\otimes_BQ^{-1}    \ar[r]^{1_N\otimes w} &   N\otimes_B W \ar[r] & 0                }
$$
Since the homomorphism  $v$ is injective, we get $s=0$. By the Snake Lemma, $t$ is surjective. Therefore  $\Coker (s)=V$ and $\Coker (t)=0$.

By Lemma \ref{sublem1}, $\mt_{B}(W) \in \Lambda_{\psi}\Gpmodcat$. Recall that $\Lambda_{\psi} \Gpmodcat$ is closed under taking kernels of surjective homomorphisms (\cite[Theorem 2.7]{holm}). Thus $(U,V,s,t) \in \Lambda_{\psi} \Gpmodcat$. By the assumption (2), we have $IU\simeq I \otimes_{\Lambda} \Coker (t)=0$ and $\Img (t)/IU \simeq N\otimes_{B}V$. Hence $\Img (t)=U\simeq N\otimes_{B}V$. This implies further that $t$ is an isomorphism by Lemma \ref{noether}(2) and that $1_N \otimes v$ is injective. Thus $\Tor_{1}^{B}(N,W)=0$.
$\square$

\medskip
(I) We show that $_BM_{\Lambda}, {}_{\Lambda}N_B$ and ${}_{\Lambda}I_{\Lambda}$ are weakly compatible.

Suppose that $\cpx{P}$ is a totally exact complex in $\mathscr{C}(\pmodcat{\Lambda})$. Then $_{\Lambda}\Ker(d_P^i)$ is Gorenstein-projective. It follows from Lemma \ref{sublem2} that $M\otimes_{\Lambda}\cpx{P}$ and $I\otimes_{\Lambda}\cpx{P}$ are exact. This further implies that the complex $\cpx{T}:=\mt_{\Lambda}(\cpx{P})=(\mt_{\Lambda}(P^i), d_T^i)$, with $d_T^i=(1_A\otimes d^i_P,1_M\otimes d^i_P))$, of $\Lambda_{\psi}$-modules is exact and $\Ker (d_T^i)=\mt_{\Lambda}(\Ker(d_P^i))$. By Lemma \ref{sublem1}, $\mt_{\Lambda}(\Ker(d_P^i))\in \Lambda_{\psi}\Gpmodcat$. Thus $\Ker (d_T^i)\in \mathscr{C}(\Gpmodcat{\Lambda_{\psi}})$ for all $i$. This implies that $\cpx{T}$ is a totally exact complex. It follows from Lemma \ref{homo-equ}(1) that $\Hom_{\Lambda}(\cpx{P},N)\simeq\Hom_{\Lambda_{\psi}}(\cpx{T},\mt_{B}(B))$ is exact.

Next, we show that $\Hom_{\Lambda}(\cpx{P},I)$ is exact. Consider the exact sequence of $\Lambda_{\psi}$-modules:
$$0\lra \mz_{\Lambda}(I) \lraf{((0,1),0)} \mt_{\Lambda} (\Lambda) \lra \mt'_{\Lambda}(\Lambda) \lra 0$$
where $\mt'_{\Lambda}(\Lambda):=(_A\Lambda, {}_BM,\mu,0)$ with $\mu: M\otimes_{A}\Lambda \ra M$ being the multiplication
 map, and $\mz_{\Lambda}(I)=(_AI, 0,0,0)\in \Lambda_{\psi}\modcat{}$. As $T^{i}:=\mt_{\Lambda}(P^i)$ is a projective $\Lambda_{\psi}$-module (see Lemma \ref{prop:projmod}(1) and $M\otimes_AP^i\simeq M\otimes_{\Lambda}P^i$), we have the exact sequence
$$0\lra \Hom_{\Lambda_{\psi}}(\cpx{T},\mz_{\Lambda}(I)) \lra  \Hom_{\Lambda_{\psi}}(\cpx{T}, \mt_{\Lambda} (\Lambda)) \lra \Hom_{\Lambda_{\psi}}(\cpx{T},\mt'_{\Lambda}(\Lambda)) \lra 0.$$
Due to the total exactness of  $\cpx{T}$, the complex $\Hom_{\Lambda_{\psi}}(\cpx{T}, \mt_{\Lambda} (\Lambda))$ is exact.  Then it follows from the exactness of $\Hom_{\Lambda}(\cpx{P},\Lambda)$ and the isomorphism $$\Hom_{\Lambda_{\psi}}(\cpx{T}, \mt'_{\Lambda}(\Lambda))=
\Hom_{\Lambda_{\psi}}(\mt_{\Lambda}(\cpx{P}), \mt'_{\Lambda}(\Lambda))\simeq \Hom_{\Lambda}(\cpx{P},\Lambda)$$ as complexes of $\mathbb{Z}$-modules that $\Hom_{\Lambda_{\psi}}(\cpx{T},\mz_{\Lambda}(I))$  is exact. This implies that $\Hom_{\Lambda}(\cpx{P},I)=\Hom_A(\cpx{P},I)$ is exact. Hence $_{\Lambda}I_{\Lambda}$ is a weakly compatible bimodule.

To complete the proof that $_BM_{\Lambda}$ and $_{\Lambda}N_B$ are weakly compatible, it remains to show that $\Hom_B(\cpx{Q},{}_BM)$ and $N\otimes_B\cpx{Q}$ are exact for any totally exact complex $\cpx{Q}\in \C{\pmodcat{B}}$. Actually, $\Tor^{1}_{B}(N,\Ker(d_Q^i))=0$ for all $i$ by Lemma \ref{sublem2}(ii). Thus $N\otimes_{B}\cpx{Q}$ is exact, and therefore $\cpx{E}:=\mt_{B}(\cpx{Q})$ is exact with $\Ker (d_{\cpx{E}}^i)=\mt_{B}(\Ker (d_Q^i))$. It follows from Lemma \ref{sublem1} that $\Ker(d_{\cpx{E}}^i)\in \Lambda_{\psi}\Gpmodcat$. Whence $\cpx{E}$ is a totally exact complex in $\mathscr{C}(\Gpmodcat{\Lambda_{\psi}})$.
Now, by Lemma \ref{homo-equ}(4), we know that $\Hom_{B}(\cpx{Q},M)\simeq \Hom_{\Lambda_{\psi}}(\cpx{E},\mt_{\Lambda}(\Lambda))$ is exact.

\medskip
(II) We prove that $\mz_{\Lambda^{op}}{(M)}$, $\mz_{\Lambda}{(N)}$, $\mz_{\Lambda^{op}}{(I)}$ and $\mz_{\Lambda}{(I)}$ are semi-weakly compatible $\Lambda_{\psi}$-modules.

Let $\cpx{T}:=(T^i,d_T^i)$ be a totally exact complex in $\mathscr{C}(\pmodcat{\Lambda_{\psi}})$. We have to show that $\mz_{\Lambda^{op}}{(M)}\otimes_{\Lambda_{\psi}}\cpx{T}$, $\Hom_{\Lambda_{\psi}}(\cpx{T},\mz_{\Lambda}{(N)})$, $\mz_{\Lambda^{op}}{(I)}\otimes_{\Lambda_{\psi}}\cpx{T}$ and
$\Hom_{\Lambda_{\psi}}(\cpx{T},\mz_{\Lambda}{(I)})$ are exact complexes.

By Lemma \ref{prop:projmod}, $\cpx{T}$ is of the form
$$  \cpx{T} \quad 	\cdots \lra \mt_{\Lambda}(P^{-1})\oplus \mt_B(Q^{-1}) \lraf{d_T^{-1}} \mt_{\Lambda}(P^{0})\oplus \mt_B(Q^{0}) \lraf{d_T^0}  \mt_{\Lambda}(P^{1})\oplus \mt_B(Q^{1})\lra \cdots.$$
As in ($\ddag$), see also the sequence (\ref{ex3}), the complex $\cpx{T}$ induces a complex $\cpx{F}=(F^i,d_F^i)\in \mathscr{C}(\pmodcat{A})$ and two complexes $\cpx{P}=(P^i,d_P^i), \cpx{Z}=(Z^i,d_Z^i)\in \mathscr{C}(\pmodcat{\Lambda})$. By Lemmas \ref{tensor}(1) and \ref{tensor}(3) and by the diagrams \ref{homo-equ}(2) and \ref{homo-equ}(7), it is sufficient to prove that $M\otimes_{\Lambda}\cpx{P}$, $\Hom_{\Lambda}(\cpx{P},N)$, $I\otimes_{\Lambda}\cpx{P}$ and $\Hom_{\Lambda}(\cpx{P},I)$ are exact complexes. Since we have shown that $_{B}M_{\Lambda}$, $_{\Lambda}N_{B}$ and $_{\Lambda}I_{\Lambda}$~ are weakly compatible, it is enough to show that $\cpx{P}$ is a totally exact complex. This is equivalent to saying that the complex $\cpx{P}$ is exact in degree $i$ and that $\Ker(d^i_P)$ is Gorenstein-projective for all $i$.

Now, we prove this statement for $i=0$.
Let $(X,Y,f,g):=\Ker (d_T^0)$, $V:=\Coker(f)$, $U:=\Coker(g)$, $(E,F,k,l):=\Ker(d_T^1)$, $T:=\Coker(k)$, $S:=\Coker(l)$.
Then $(X,Y,f,g)$ and $(E,F,k,l)$ are Gorenstein-projective $\Lambda_{\psi}$-modules. Consider the canonical exact sequence  of $\Lambda_{\psi}$-modules:
$$0\lra (X,Y,f,g) \stackrel{(d_X,d_Y)}{\hookrightarrow} T^{0} \lra (E,F,k,l) \lra 0$$
with  $d_X: {}_AX\ra {}_AP^0(I)\oplus {}_AN\otimes_{B}Q^0$, $d_Y: {}_BY\ra {}_BM\otimes_{\Lambda}P^0\oplus Q^0$. Then $(d_X,d_Y)$ is a left $\add({}_{\Lambda_{\psi}}\Lambda_{\psi})$-approximation of $(X.Y,f,g)$. By the diagram (\ref{morphismdiagrams}), there is the exact commutative diagram of $A$-modules
\begin{equation}\label{ex13}
\xymatrix{
  	& N\otimes_B Y \ar[r]\ar[d]^{g}& N\otimes_BM\otimes_{\Lambda}P^{0} \oplus N\otimes_B Q^0 \ar[r] \ar[d]^{\left(\begin{smallmatrix}
   \psi_{P^0} & 0 \\
   0 & 1_{N\otimes_B Q^0 }
\end{smallmatrix}\right)} & N\otimes_B F \ar[r]\ar[d]^{l} & 0\\
  	0\ar[r]& X \ar[r]^-{d_X}\ar[d]^{\lambda_X}& P^{0}(I) \oplus N\otimes_{B}Q^{0} \ar[r] \ar[d] & E \ar[r]\ar[d]& 0 \\
  	& U\ar[r]^{a_U}\ar[d] & P^{0} \ar[r]^{b_S}\ar[d]& S\ar[r]\ar[d] & 0\\
  	& 0 &0 &0}
\end{equation}
Note that the bottom row is exact by the Snake Lemma, and thus an exact sequence of $\Lambda$-modules by Lemma \ref{module-A}(1). Let $d_{{}_{\Lambda}X}=(e_0,e_1,e_2):{}_{\Lambda}X \ra P^0 \oplus I\otimes_{\Lambda}P^0 \oplus N\otimes Q^0$ denote the restriction of $d_X$ to $\Lambda$-modules. Then the above diagram shows $e_0=\lambda_Xa_U$.

Similarly,  let $\Ker(d^{-1}_T) =(E',F', k',l')$ and $(t_X,t_Y): T^{-1}\to \Ker(d_T^0)=(X,Y,f,g)$ be the canonical projection. Then there is a canonical exact sequence of $\Lambda_{\psi}$-modules
$0\lra (E',F',k',l') \lra T^{-1} \lraf{(t_X,t_Y)} (X,Y,f,g) \lra 0.$ This supplies us with the exact commutative diagram of $A$-modules
\begin{equation}\label{ex14}
\xymatrix{
  	& N\otimes_B F' \ar[r]\ar[d]^{l'}& N\otimes_BM\otimes_{\Lambda}P^{-1} \oplus N\otimes_B Q^{-1} \ar[r] \ar[d]^{\left(\begin{smallmatrix}
   \psi_{P^{-1}} & 0 \\
   0 & 1_{N\otimes_B Q^{-1} }
\end{smallmatrix}\right)} & N\otimes_B Y \ar[r]\ar[d]^{g} & 0\\
  	0\ar[r]& E' \ar[r]\ar[d]& P^{-1}(I) \oplus N\otimes_{B}Q^{-1} \ar[r]^-{t_X} \ar[d] & X \ar[r]\ar[d]^{\lambda_X}& 0 \\
  	& S'\ar[r]^{a_{S'}} & P^{-1} \ar[r]^{b_U}& U\ar[r] & 0\\
}
 \end{equation}
Due to $d_T^{-1}=(t_X,t_Y)(d_X,d_Y)$, we have $t_Xd_X=d_F^{-1}$. Thus the diagrams (\ref{ex13}) and (\ref{ex14}) provide the commutative diagram of $A$-modules:
   $$ \xymatrix@C=0.6cm@R=0.6cm{P^{-1}(I) \oplus N\otimes_{B}Q^{-1}
 \ar[r]^-{t_X}\ar[d] & X   \ar[r]^-{d_X}\ar[d]^{\lambda_{X}} &  P^{0}(I)\oplus N\otimes_{B}Q^{0}\ar[d] \\
   P^{-1}  \ar@{->>}[r]^{b_U} & U  \ar[r]^{a_U} & P^{0}}
 $$
where the two unnamed vertical maps are natural projections. Since $d_F^{-1}= t_Xd_X$, we get $b_Ua_U=d_P^{-1}$.

We show that $a_U$ is injective. Consider $a_U$ as a homomorphism of $\Lambda$-modules. Since $(X,Y, f, g)$ is Gorenstein-projective, we know $U\in \Lambda\Gpmodcat$ by the assumption (2). Thus, for $a_U$ to be injective, it suffices to show that $a_U$ is a left $\add({}_{\Lambda}\Lambda)$-approximation of $U$. Actually, for $P\in \add({}_{\Lambda}\Lambda)$ and a homomorphism $a_0:U\ra P$ of $\Lambda$-module, we have to find a homomorphism $e: P^0\to P$ of $\Lambda$-modules, such that $a_0=a_Ue$. To define $e$, we will construct
a $\Lambda$-module homomorphism $a_1:{}_{\Lambda}X\ra I\otimes_{\Lambda}P$, where $X$ is regarded as $\Lambda$-module by restriction, and a $B$-module homomorphism $h: Y\ra M\otimes_{\Lambda}P$, such that $a:=(\lambda_{X}a_0,a_1)$ is a homomorphism of $A$-modules, and  $(a,h): (X, Y,f,g)\ra \mt_{\Lambda}(P)$ is a homomorphism of $\Lambda_{\psi}$-modules.

\textbf{Step 1.} Construction of $h$. By Lemma \ref{module-A}(2), we have a commutative diagram
$$\xymatrix@C=0.6cm@R=0.5cm{
_BM\otimes_{A}X \ar[r]^{\lambda'_X} \ar[d]^{f} & _BM\otimes_{\Lambda}U \ar[dl]^{\eta_Y}\ar[r]& 0\\
   _BY}$$
where $\lambda'_X$ is the composite of $1_M\otimes\lambda_X$ with the isomorphism $M\otimes_{A}U\ra M\otimes_{\Lambda}U$ as $B$-modules. By assumption, $\Img(f)\simeq M\otimes_{\Lambda}U$. This means that $\eta_Y$ is injective and there is an exact sequence of $B$-modules
$$   0\lra M\otimes_{\Lambda}U \lraf{\eta_Y} Y \lraf{\mu_Y} V \lra 0. $$
Since ${}_{B}V \in B \Gpmodcat$ and ${}_{B}M$ is semi-weakly compatible, there holds $\Ext^{1}_{B}(V,M\otimes_{\Lambda}P)=0$. This shows that $\Hom_{B}(\eta_{Y},M\otimes_{\Lambda}P)$ is surjective, and therefore there is a homomorphism $h: Y \ra M\otimes_{\Lambda}P$ of $B$-modules such that $\eta_{Y}h=1_M\otimes a_0$.

\textbf{Step 2.} Construction of $a_1$. From $\eta_{Y}h=1_M\otimes a_0$, one gets $(1_N\otimes_B\eta_Y)(1_N\otimes h)=1_N\otimes_B 1_M \otimes a_0$. It follows from the natural property of $\psi$  that the diagram of $A$-modules is commutative
$$\xymatrix{
N\otimes_{B}M\otimes_{\Lambda}U \ar[r]^-{1_N\otimes \eta_Y} \ar[d]^{\psi\otimes 1_U} & N\otimes_{B}Y \ar[d]^{(1_N\otimes h)(\psi\otimes 1_P)}\\
I\otimes_{\Lambda}U  \ar[r]^{1_I \otimes a_0} & I\otimes_{\Lambda}P
}$$

Now, let $H:=\Img(g)$ and $g = \sigma \epsilon_X$ with $\sigma: N\otimes_BY\ra H$ the canonical projection and $\epsilon_X: H\hookrightarrow X$ the inclusion. According to Lemma \ref{b-condition} and its proof, there exists an injective homomorphism $m: {}_AI\otimes_{\Lambda}U\to {}_AH$, such that the following is a pushout diagram
   $$
\xymatrix{
&N\otimes_{B}M\otimes_{\Lambda}U \ar[r]^-{1_N\otimes \eta_Y} \ar[d]^{\psi\otimes 1_U} & N\otimes_{B} Y \ar[d]^{\sigma} \\
0\ar[r] & I\otimes_{\Lambda}U \ar[r]^{m} & H  }
$$
By a universal property of pushouts, there is a $\Lambda$-module homomorphism $t: H \ra I\otimes_{\Lambda}P$, such that $1_I\otimes a_0 = mt$ and $(1_N\otimes h)(\psi\otimes 1_P) = \sigma\, t$.

The exact sequence $ 0\ra H\lraf{\epsilon_X} X\lraf{\lambda_X} {}_AU\ra 0$ of $A$-modules restricts to an exact sequence of $\Lambda$-modules:
$$ 0 \lra   H \lraf{\epsilon_X}{}_{\Lambda}X \lraf{\lambda_X} U \lra 0 $$
It follows from $_{\Lambda}U \in \Lambda \Gpmodcat$ and the semi-weak compatibility of $_{\Lambda}I$ that $\Ext^{1}_{\Lambda}(U,I\otimes_{\Lambda}P)=0$, and therefore $(\epsilon_{X},I\otimes_{\Lambda}P)$ is surjective. Hence there is a homomorphism $a_1: X \ra I\otimes_{\Lambda}P$ of $\Lambda$-modules such that $\epsilon_{X}a_1=t$.

\textbf{Step 3.} We show that $(\lambda_Xa_0,a_1): X\ra P\oplus I\otimes_{\Lambda}P$ is a homomorphism of $A$-modules. We write $a$ for $(\lambda_Xa_0,a_1)$ for simplicity. On the one hand, $1_I\otimes a_0=mt=m\epsilon_{X}a_1$.
On the other hand, by the definition of $m$ (see Lemma \ref{b-condition}), the following diagram commutes
$$\xymatrix@C=0.8cm@R=0.5cm{I\otimes_{A}X \ar[r]^-{\mlt_X} \ar[d]_-{1_I\otimes\lambda_X}  &X\\
    I\otimes_{\Lambda}U \ar[r]^{m} & H \ar[u]_{\epsilon_X}}$$
that is, $\mlt_X= (1_I\otimes \lambda_X)\, m \, \epsilon_X$. Thus, for $i\in I$ and $x\in X$, there holds
$$(ix)a_1= \big((i\otimes x)\mlt_X\big)a_1=[(i\otimes x)\big((1_I\otimes\lambda_X)\, m\, \epsilon_X\big)]a_1=(i\otimes x)[(1_I\otimes\lambda_X)(1_I\otimes a_0)] = i\otimes (x)\lambda_Xa_0.$$
Clearly, $(ix)\lambda_Xa_0=0$. Now it is easy to verify that $a$ is a homomorphism of $A$-modules.

\textbf{Step 4.} We prove that $(a,h)$ is a homomorphism of $\Lambda_{\psi}$-modules. First, we show that the out-side square of the following diagram of $A$-modules is commutative
   $$\xymatrix@C=0.8cm@R=0.6cm{
   N\otimes_{B}Y \ar[rr]^{1_N\otimes h} \ar[dd]^{g}\ar[dr]^-{\sigma} & &N\otimes_{B}M\otimes_{\Lambda}P \ar[dd]^{\psi_P}\\
   & H \ar[dl]^{\epsilon_X}\ar[dr]^-{(0,t)}\\
   X \ar[rr]^{a} & &P(I)
   }$$
In fact, since $\epsilon_X\lambda_Xa_0=0$, we have $\epsilon_Xa=(0,\epsilon_Xa_1)=(0,t)$ by Step 2. It follows from $(1_N \otimes h)(\psi\otimes 1_P)=\sigma\, t$ that $(1_N \otimes h)\psi_P
=\sigma(0,t)$. Then $ga=(\sigma \epsilon_X) a=\sigma (0,t) =(1_N\otimes h)\psi_P$.

Second, we show that the out-side square of the following diagram of $B$-modules is commutative
$$\xymatrix@C=0.8cm@R=0.6cm{
   M\otimes_{A}X \ar[rr]^{1_M\otimes a} \ar[dd]^{f}\ar[dr]^-{\lambda'_X} & & M\otimes_{A}P(I) \ar[dd]^{\pi_P}\\
   &  M\otimes_{\Lambda}U \ar[dl]^{\eta_Y}\ar[dr]^-{1_M \otimes a_0}\\
   Y \ar[rr]^{h} & &  M\otimes_{\Lambda}P
   }$$
Note that $\lambda'_X\eta_Y=f$ and $\eta_Yh=1_M\otimes a_0$. A straightforward verification shows $fh=(1_M\otimes a)\pi_P$. Thus the pair $(a,h)$ is a homomorphism of $\Lambda_{\psi}$-modules.

\textbf{Step 5.} Definition of $e$.
Since $(d_X,d_Y)$ is a left $\add(_{\Lambda_{\psi}}\Lambda_{\psi})$-approximation of $(X,Y,f,g)$ and $\mt_{\Lambda}(P)\in \add(_{\Lambda_{\psi}}\Lambda_{\psi})$, there is a homomorphism $u:\mt_{\Lambda}(P^0) \ra \mt_{\Lambda}(P)$ and a homomorphism $v:
\mt_{B}(Q^0) \ra \mt_{\Lambda}(P)$ such that $(d_X,d_Y)\binom{u}{v}=(a,h)$. By Lemmas \ref{homo-equ}(3)-(4), there exists a homomorphism $e:P\ra P^0$ of $\Lambda$-modules, satisfying $e_0e=\lambda_Xa_0$. Because $e_0=\lambda_Xa_U$ and $\lambda_X$ is surjective, there holds $a_Ue=a_0$. Hence $a_U$ is a left $\add({}_{\Lambda}\Lambda)$-approximation of $_{\Lambda}U$. This completes the proof of $a_U$ being injective.

\medskip
Now we show that the complex $\cpx{P}$ is totally exact. From the exact sequence
$$0\lra U \lraf{a_U}  P^0 \lraf{b_S} S \lra 0$$of $\Lambda$-modules, we proceed with a similar proof of $a_U$ being a left $\add(_\Lambda\Lambda)$-approximation of $U$ with $d_P^{-1}=b_Ua_U$, and replace $U$ with $S$  to show that there is an injective homomorphism $a_S: S\to P^1$ such that $d_P^0=b_Sa_S$. This implies that $\Ker (d_P^0)=\Ker (b_S)=\Img(a_U)\simeq U\in \Lambda \Gpmodcat$. Due to $\Img(d_P^{-1})= \Img(a_U)$, we see that $\Ker(d^0_P)=\Img(d^{-1}_P)$ and $\cpx{P}$ is exact in degree $0$ with $\Ker(d_P^0)\in  \Lambda\Gpmodcat$. Similarly, we can show that $\cpx{P}$ is exact in any degree $i$ with $\Ker (d_P^i)\in \Lambda \Gpmodcat$. Thus $\cpx{P}$ is a totally exact complex. This finishes the proof of (2) implying (1). $\square$

\medskip
For the special Morita context ring $\Lambda_{(0,0)}$, it was shown in \cite{GaoP17} that the compatibility conditions suffice  a class of modules over $\Lambda_{(0,0)}$ to be Gorenstein-projective. Next, we shall point out that the weak compatibility conditions are both sufficient and necessary.

\begin{Prop} \label{zero-case} For the Morita context ring $\Lambda_{(0,0)}$, the following are equivalent.

$(1)$ $_AN_B$ and $_BM_A$ are weakly compatible bimodules; $(_AN,0,0,0)$ and $(M_A, 0,0,0)$ are semi-weakly compatible left and right $\Lambda_{(0,0)}$-modules, respectively.

$(2)$ A $\Lambda_{(0,0)}$-module $(X,Y,f,g)$ is Gorensteion-projective if and only if

$\quad (a)$  $_B\Coker(f)$ and $_A\Coker(g)$ are Gorenstein-projective, and	

$\quad (b)$	$_B\Img(f)\simeq {}_BM\otimes_A\Coker(g)$ and $_A\Img(g)\simeq {}_AN\otimes_B\Coker(f)$,where $\Coker(f)$ and $\Img(g)$ denote the cokernel of $f$ and the image of $g$, respectively.
\end{Prop}

{\bf Proof.} This follows immediately from Theorem \ref{main-result-3} because $I=0$ in $\Lambda_{(0,0)}$. $\square$

\medskip
Remark that the following was proved in \cite[Theorem A(i)]{GaoP17}: Assume that both $M$ and $N$ are compatible bimodules over Artin algebras. If a $\Lambda_{(0,0)}$-module $(X,Y,f,g)$ fulfills the conditions $(a)$ and $(b)$ in Proposition \ref{zero-case}, then $(X,Y,f,g)$ is Gorensteion-projective. It seems that the weak compatibility conditions are more suitable for describing Gorenstein-projective modules over $\Lambda_{(0,0)}$.

\section{Applications to noncommutative tensor products\label{sect4}}
In this section we will construct Gorenstein-projective modules over noncommutative tensor products of exact contexts arising from Morita contexts with two bimodule homomorphisms zero. This is related to the Morita context rings $\Lambda_{(\phi,0)}$.

\begin{Def} {\rm \cite{xc3}} Let $\lambda: R\to S, \mu: R\to T$ be homomorphisms of unitary rings, and $_SW_T$ an $S$-$T$-bimodule with $w\in W$. If the sequence
$$ 0\lra R\lraf{(\lambda,\mu)} S\oplus T\lraf{\binom{\cdot w}{w\cdot}} W\lra 0$$ is exact of abelian groups, then $(\lambda, \mu, W, w)$ is called an exact context, where $\cdot w: S\to W$ is the right multiplication by $w$. The noncommutative tensor product of $(\lambda, \mu, W, w)$ is well defined.
\end{Def}

Morita contexts provide prominent examples of exact contexts. For a Morita context $(A, \Gamma, {}_{\Gamma}M_A, {}_AN_{\Gamma}, \phi, \psi)$, let
$$ R:=\left(\begin{array}{lc} A &0 \\
0& \Gamma \end{array}\right),\; S:=\left(\begin{array}{lc} A &N \\
0& \Gamma\end{array}\right),\; T:=\left(\begin{array}{lc} A &0 \\
M& \Gamma\end{array}\right), \; W:=\Lambda_{(\phi,\psi)}, \;  w := \binom{1 \quad 0}{0 \quad 1}.$$
If $\lambda$ and $\mu$ are the inclusions, then $(\lambda, \mu, W, w)$ is an exact context. Its noncommutative tensor product, denoted by $C(A,\Gamma, M, N,\phi,\psi)$, can be described explicitly: $C(A, \Gamma, M, N,\phi,\psi)$ has the underlying abelian group of the matrix form
$$\left(\begin{array}{lc} A & N \\
M & \Gamma\oplus (M\otimes_AN)\end{array}\right)$$ with the multiplication $\circ$ defined by
$$\left(\begin{array}{lc} a_1 & n_1\\
m_1& (b_1, m\otimes n)\end{array}\right)\circ \left(\begin{array}{lc} a_2 &n_2 \\
m_2& (b_2, m'\otimes n')\end{array}\right)$$ $$ = \left(\begin{array}{lc} a_1a_2+(n_1\otimes m_2)\psi &a_1n_2+n_1b_2+n_1(m'\otimes n')\phi \\
m_1a_2+b_1m_2+(m\otimes n)\phi m_2& \big(b_1b_2,\,m_1\otimes n_2+(b_1m')\otimes n'+m\otimes (nb_2)+ m\otimes (n\otimes m')\psi n'\big)\end{array}\right),$$
where $a_1, a_2\in A, b_1, b_2\in \Gamma, n_1, n_2, n, n'\in N$ and $m_1, m_2, m, m'\in M$. For details, we refer the reader to \cite{xc3}.

Let $C:=C(A,\Gamma, M, N,0,0)$, and let $B:=\Gamma \ltimes (M\otimes_AN)$ be the trivial extension of $\Gamma$ with the $\Gamma$-bimodule $M\otimes_A N$. We may regard $M$ as a $B$-$A$-bimodule and $N$ as an $A$-$B$-bimodule via the canonical surjective homomorphism $B\to \Gamma$. Thus we have a Morita context $(A,B,M,N,\phi, 0)$, where $\phi: M\otimes_AN\to B, m\otimes n\mapsto (0, m\otimes n)$ for $m\in M, n\in N$, and the Morita context ring $\Lambda_{(\phi,0)}=\begin{pmatrix} A & N\\ M & B\end{pmatrix}_{(\phi,0)}$ which is isomorphic to $C$. Thus the dual version of Theorems \ref{main-result-1} and \ref{main-result-2} describe also the Gorenstein-projective modules over the noncommutative tensor product $C$. For example, we have the specifical corollary.

\begin{Koro} Suppose that the bimodule $_{\Gamma}M_A, {}_AN_{\Gamma}$ and $ {}_{\Gamma}M\otimes_AN_{\Gamma}$ are weakly compatible. Let $B=\Gamma\ltimes J$  with $J:=M\otimes_AN$. Then a $C$-module $(_AX,{}_BY,f,g)$ is Gorenstein-projective if

$(i)$ ${}_{\Gamma}\Coker(f)$ and $_A\Coker(g)$ are Gorenstein-projective, and

$(ii)$ $_AN\otimes_B\Coker(f)\simeq {}_A\Img(g)$, $_BM\otimes_A\Coker(g)\simeq {}_B\Img(f)/JY$, and
$_BJ\otimes_B\Coker(f)\simeq {}_BJY$, where $\Coker(f)$ and $\Img(g)$ denote the cokernel of $f$ and the image of $g$, respectively.
\label{nonctp}
\end{Koro}

\medskip
{\bf Acknowledgements.} The research work was supported partially by the National Natural Science Foundation of China (Grant 12031014).

\medskip
{\footnotesize

Qianqian Guo, School of Mathematical Sciences, Capital Normal University, 100048 Beijing, P.R. China

{\tt Email: 1559147300@qq.com (Q.Q. Guo)}

\medskip
Changchang Xi, School of Mathematical Sciences, Capital Normal University, 100048 Beijing, P.R. China; and
School of Mathematics and Information Science, Henan Normal University, 453007 Xinxiang, Henan, P. R. China

{\tt Email: xicc@cnu.edu.cn (C.C. Xi)}
}
\end{document}